\documentclass{amsart}
\usepackage[dvips]{graphicx}
\usepackage{amssymb,amsmath}
\usepackage[all]{xy}

\allowdisplaybreaks[1]

\newtheorem{theorem}{Theorem}[section]
\newtheorem{proposition}[theorem]{Proposition}
\newtheorem{lemma}[theorem]{Lemma}

\theoremstyle{definition}
\newtheorem{definition}[theorem]{Definition}

\theoremstyle{remark}
\newtheorem{remark}[theorem]{Remark}

\numberwithin{equation}{section}

\newcommand{\PSLC}{\mathrm{PSL}(2,\mathbb{C})}

\newcommand{\Log}{\mathrm{Log}}

\newcommand{\hatPC}{\widehat{\mathcal{P}}(\mathbb{C})}
\newcommand{\hatBC}{\widehat{\mathcal{B}}(\mathbb{C})}
\newcommand{\hatP}{\widehat{\mathcal{P}}}

\newcommand{\hatbeta}{\widehat{\beta}}
\newcommand{\Cpm}{(\mathbb{C}^{2} \setminus \{ 0 \}) / \pm}

\begin{document}

\title[Quandle homology and complex volume]{Quandle homology and complex volume}

\author[A. Inoue and Y. Kabaya]{Ayumu Inoue and Yuichi Kabaya}
\address{
Department of Mathematics Education, Aichi University of Education,
Kariya, Aichi 448--8542 Japan}
\email{ainoue@auecc.aichi-edu.ac.jp} 
\address{
Department of Mathematics, Osaka University, 
Toyonaka, Osaka 560--0043 Japan}
\email{y-kabaya@cr.math.sci.osaka-u.ac.jp}


\begin{abstract}
We introduce a new homology theory of quandles, called simplicial quandle homology, which is quite different from quandle homology developed by Carter et al.
We construct a homomorphism from a quandle homology group to a simplicial quandle homology group.
As an application, we obtain a method for computing the complex volume of a hyperbolic link only from its diagram.

\end{abstract}

\maketitle

\section{Introduction}
\label{sec_intro}
The volume and the Chern-Simons invariant of a finite volume hyperbolic 3-manifold have been intensively studied.
In this paper we give a description of these invariants for link complements from the view point of quandle cocycles.
This enables us to compute these invariants in terms of a link diagram and its coloring.

The Cheeger-Chern-Simons class $\hat{c}_2$ is an invariant of flat $\PSLC$-bundles.
This gives a homomorphism from the group homology $H_3(\PSLC; \mathbb{Z})$ as a discrete group to $\mathbb{C} / \pi^2 \mathbb{Z}$.
For a closed oriented 3-manifold and a representation of its fundamental group into $\PSLC$, 
we obtain an invariant with values in  $\mathbb{C}^2 / \pi^2 \mathbb{Z}$ by evaluating at the fundamental class of the manifold.
For a finite volume hyperbolic 3-manifold, there is a unique discrete faithful representation up to conjugation, 
in this case $\hat{c}_2$ is equal to $i(\mathrm{Vol}+i\mathrm{CS})$, where $\mathrm{Vol}$ and $\mathrm{CS}$ are volume and Chern-Simons invariant of the hyperbolic metric.

Dupont gave a formula of the Cheeger-Chern-Simons class modulo $\pi^2 \mathbb{Q}$ by using the Rogers dilogarithm function in \cite{dupont}.
Then Neumann gave a formula modulo $\pi^2$ in \cite{neumann04}.
He defined the extended Bloch group $\hatBC$ and showed that $\hatBC$ is isomorphic to the group homology $H_3(\PSLC;\mathbb{Z})$.
He also defined a map $R: \hatBC \to \mathbb{C}/\pi^2\mathbb{Z}$ which gives the Cheeger-Chern-Simons class.
To apply his formula to a hyperbolic $3$-manifold, we need an ideal triangulation of the $3$-manifold. 
Quandle homology plays a useful role to give a diagrammatic description of these invariants.

A quandle, which was introduced by Joyce in \cite{joyce}, is a set with a binary operation like conjugation in a group.
For a quandle $X$, Carter et al. defined a quandle homology $H^{Q}_{n}(X)$ in \cite{CJKLS03}.
Then various kinds of generalization has been introduced by several authors.
In this paper we introduce a new quandle homology $H^{\Delta}_n(X)$, which we call simplicial quandle homology.
We shall construct a map from a quandle homology $H^Q_n(X;\mathbb{Z}[X])$ to the simplicial quandle homology $H^{\Delta}_{n+1}(X)$.
Roughly speaking, this map gives a triangulation of a link complement.

For a diagram of an oriented link $L$ in $S^3$, we can define the notion of shadow colorings by a quandle $X$.
A shadow coloring is a pair of maps $\mathcal{S}=(\mathcal{A},\mathcal{R})$, where $\mathcal{A}$ is a map from the set of arcs of the diagram to $X$
and $\mathcal{R}$ is a map from the set of complementary regions to $X$ satisfying some conditions.
For a shadow coloring $\mathcal{S}$, we define a cycle $[C(\mathcal{S})]$ of $H^Q_2(X; \mathbb{Z}[X])$.
Let $\mathcal{P}$ be the set of parabolic elements of $\PSLC$, then $\mathcal{P}$ has a quandle structure by conjugations.
There exists a one-to-one correspondence between the set of arc colorings by $\mathcal{P}$ 
and the set of $\PSLC$-representations of the fundamental group $\pi_{1}(S^{3} \setminus L)$ which sends each meridian to a parabolic element.
We will show that the homology class $[C(\mathcal{S})]$ with respect to $\mathcal{P}$ does only depend on the conjugacy class of the representation.

We will see that the quandle $\mathcal{P}$ is identified with $\Cpm$, the set of non-zero two-dimensional complex vectors identifying $v$ with $-v$.
From this fact, $H^{\Delta}_3(\mathcal{P})$ is closely related to the homology group $H_3(C^{h \neq}_*( \mathbb{C}^2)_{\PSLC})$, 
which was studied by Dupont and Zickert in \cite{dupont-zickert}.
Following their construction we obtain a homomorphism
\begin{equation}
\label{eq:intro_homo}
 H^Q_2(\mathcal{P}; \mathbb{Z}[\mathcal{P}]) \longrightarrow  \hatBC.
\end{equation}
We will show that the image of $[C(\mathcal{S})]$ under the map gives the invariant defined by Neumann (Theorem \ref{thm:volume_and_chern_simons}).
Composing with the map $R : \hatBC \to \mathbb{C} / \pi^2 \mathbb{Z}$, we obtain 
\[
\mathrm{cvol} :H^Q_2(\mathcal{P}; \mathbb{Z}[\mathcal{P}]) \longrightarrow  \mathbb{C} / \pi^2 \mathbb{Z}.
\]
This gives a cocycle $[\mathrm{cvol}] \in H_Q^2(\mathcal{P}; \mathrm{Hom}(\mathcal{P}, \mathbb{C} / \pi^2 \mathbb{Z}) )$.
For a hyperbolic link $L$ and a shadow coloring $\mathcal{S}$ which corresponds to the discrete faithful representation, we have 
\begin{equation}
\label{eq:intro_main}
\langle [\mathrm{cvol}], [C(\mathcal{S})]\rangle = i (\mathrm{Vol}(S^3 \setminus L) + i \mathrm{CS}(S^3 \setminus L)),
\end{equation}
where $\mathrm{CS}(S^{3} \setminus L)$ is the Chern-Simons invariant for cusped hyperbolic manifolds defined by Meyerhoff \cite{meyerhoff}.

Since $[C(\mathcal{S})]$ is an invariant of conjugacy classes of representations,
$\langle [\mathrm{cvol}], [C(\mathcal{S})] \rangle$ is clearly a link invariant with values in $\mathbb{C} / \pi^2 \mathbb{Z}$ 
when the coloring corresponding to the discrete faithful representation,  
even if the reader does not know the definition of the hyperbolic volume and the Chern-Simons invariant.
Our description is based on the quandle homology theory, but we do not care about it in actual calculations
as explained in Section \ref{sec_example}.

Finally, we remark that the triangulation of a link complement constructed in this paper is not an ideal triangulation in the usual sense.
It may contain negatively oriented or flat ideal tetrahedra, and non-ideal vertices.
See Remark \ref{rem:on_the_tirangulation} for details.

This paper is organized as follows.
In Section \ref{sec_quandle_(co)homology}, we recall the definition of quandles and quandle homology theory.
We introduce the simplicial quandle homology in Section \ref{sec_main} and construct a homomorphism 
from a usual quandle homology to the simplicial quandle homology.
In Section \ref{sec_shadow_coloring}, we review the definition of a shadow coloring and the cycle associated with a shadow coloring.
We study the quandle $\mathcal{P}$ consisting of parabolic elements of $\PSLC$ in Section \ref{sec_psl2c}.
In Section \ref{sec_extended_bloch}, we recall the definition of the extended Bloch group defined by Neumann.
Then we construct the homomorphism (\ref{eq:intro_homo}) in Section \ref{sec_invariant}.
In the final section, we demonstrate a computation of the volume and the Chern-Simons invariant for the $5_2$ knot.


\subsection*{Acknowledgments}
The authors would like to express their sincere gratitude to Professor Sadayoshi Kojima for encouraging them.
The first author was supported in part by JSPS Global COE program ``Computationism as a Foundation for the Sciences''.
The second author was partially supported by JSPS Research Fellowships for Young Scientists.
We thank Christian Zickert for for useful conversations. Finally, we also thank 
the referees for helpful comments.

\section{Quandle and quandle homology}
\label{sec_quandle_(co)homology}
In this section, we recall the definitions of a quandle and quandle homology.

\subsection{Quandle}
A \emph{quandle} is a non-empty set $X$ equipped with a binary operation $\ast$ satisfying the following three axioms:
\begin{itemize}
 \item[(Q1)] For any $x \in X$, $x \ast x = x$.
 \item[(Q2)] For any $y \in X$, the map $\ast y : X \rightarrow X$ ($x \mapsto x \ast y$) is bijective.
 \item[(Q3)] For any $x, y, z \in X$, $(x \ast y) \ast z = (x \ast z) \ast (y \ast z)$.
\end{itemize}
For example, let $X$ be a subset of a group closed under conjugations.
Then $X$ is a quandle with $x \ast y = y^{-1} x y$ for any $x, y \in X$.
We call it a \emph{conjugation quandle}.
The notions of homomorphisms and isomorphisms of quandles are appropriately defined.

Let $X$ be a quandle.
Define a binary operation $\ast^{-1}$ of $X$ so that the map $\ast^{-1} y : X \rightarrow X$ is the inverse of the bijection $\ast y : X \rightarrow X$ for any $y \in X$.
Then $\ast^{-1}$ also satisfies the three axioms of a quandle.

The \emph{associated group} $G_{X}$ of $X$ is a group generated by elements $x \in X$ subject to the relation $x \ast y = y^{-1} x y$ for each $x, y \in X$.

Suppose $g = x_{1}^{\varepsilon_{1}} x_{2}^{\varepsilon_{2}} \cdots x_{n}^{\varepsilon_{n}}$ is an element of $G_{X}$ 
with some $n \geq 0$, $x_{i} \in X$, and $\varepsilon_{i} \in \{ \pm 1 \}$.
For each $x \in X$, define an element $x \ast g \in X$ by
\[
 x \ast g = ( \cdots ((x \ast^{\varepsilon_{1}} x_{1}) \ast^{\varepsilon_{2}} x_{2}) \cdots ) \ast^{\varepsilon_{n}} x_{n}.
\]
Here, $\ast^{+1}$ denotes the binary operation $\ast$.
Then a map $X \times G_{X} \rightarrow X$ sending $(x, g)$ to $x \ast g$ is a right action of $G_{X}$ on $X$.

\subsection{Quandle homology}
\begin{figure}[htb]
\begin{center}
\includegraphics[width=110pt,clip]{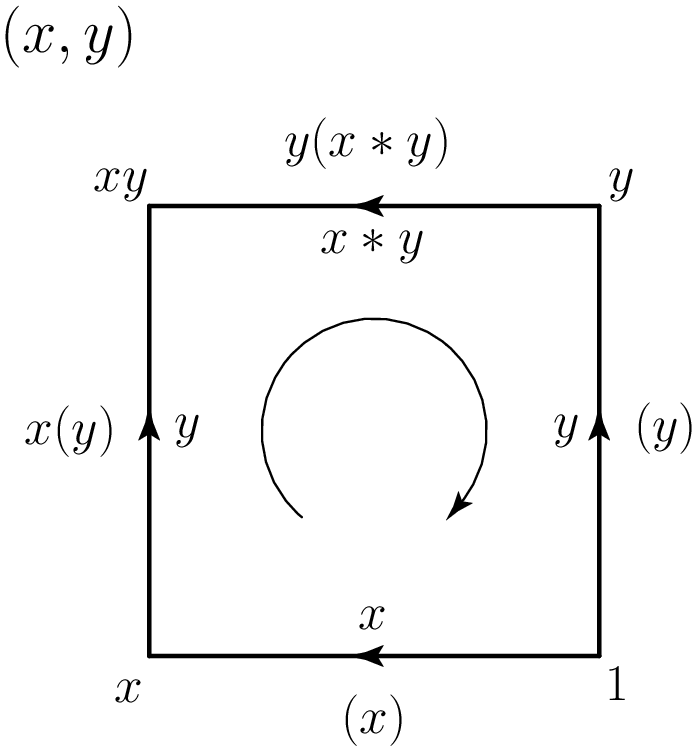}
\hspace{20pt}
\includegraphics[width=170pt,clip]{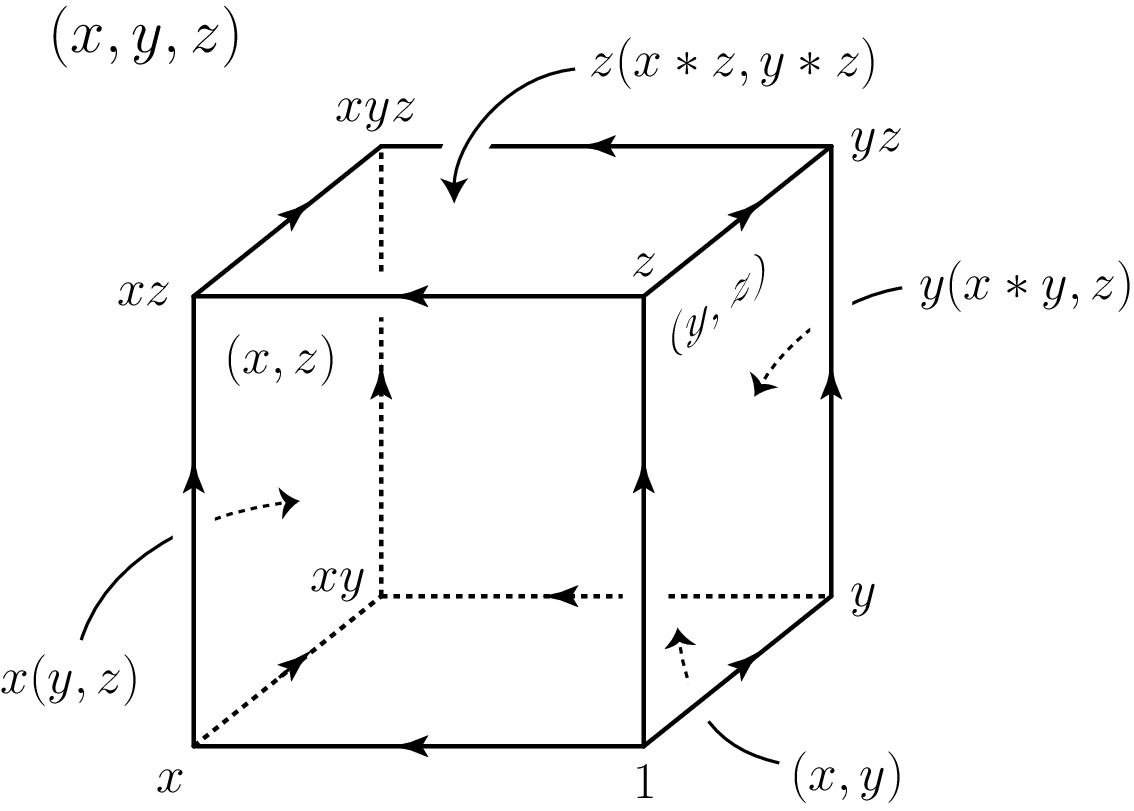}
\end{center}
\caption{Schematic picture of the boundary map:
The left indicates $\partial (x,y) = -(y) + x(y) +(x) -y(x*y)$. 
The right indicates $\partial (x,y,z) = -(y,z) +x(y,z) +(x,z) -y(x*y,z) -(x,y) +z(x*z,y*z)$. 
We use a non-standard orientation convention to be consitent with the positive crossing in Fig. \ref{fig:cycle_from_shadow_coloring}.
}
\label{fig:boundary_map}
\end{figure}

Let $X$ be a quandle, $G_{X}$ the associated group of $X$, and $\mathbb{Z}[G_{X}]$ the group ring of $G_{X}$.
Consider the free left $\mathbb{Z}[G_{X}]$-module $C^{R}_{n}(X)$ generated by all $n$-tuples $(x_{1}, x_{2}, \cdots, x_{n}) \in X^{n}$ for each $n \geq 1$.
We let $C^{R}_{0}(X) = \mathbb{Z}[G_{X}]$.
Define a map $\partial : C^{R}_{n}(X) \rightarrow C^{R}_{n-1}(X)$ by
\begin{eqnarray*}
 \partial(x_{1}, x_{2}, \cdots, x_{n}) = \sum_{i = 1}^{n} (-1)^{i} \{(x_{1}, x_{2}, \cdots, \widehat{x_{i}}, \cdots, x_{n}) & & \\
 & & \hskip -12.5em - \> x_{i} (x_{1} \ast x_{i}, x_{2} \ast x_{i}, \cdots, x_{i-1} \ast x_{i}, x_{i+1}, \cdots, x_{n}) \}.
\end{eqnarray*}
Then $\partial$ satisfies $\partial \circ \partial = 0$.
Therefore, $C^{R}_{\ast}(X) = (C^{R}_{n}(X), \partial)$ is a chain complex.
As illustrated in Figure \ref{fig:boundary_map}, a generator of $C^{R}_{n}(X)$ may be identified with an $n$-cube whose edges are labelled by elements of $X$ and vertices are labelled by elements of $G_{X}$.
The boundary map $\partial$ sends a cube to a formal sum of its boundaries.

Define a submodule $C^{D}_{n}(X)$ of $C^{R}_{n}(X)$ by
\[
 C^{D}_{n}(X) = 
 \begin{cases}
  \mathrm{span}_{\mathbb{Z}[G_{X}]} \{ (x_{1}, x_{2}, \cdots, x_{n}) \in X^{n} \mid x_{i} = x_{i+1} \ \mathrm{for \ some} \ i \} & (n \geq 2), \\
  0 & (n = 0, 1).
 \end{cases}
\]
It is routine to check that $C^{D}_{\ast}(X) = (C^{D}_{n}(X), \partial)$ is a subchain complex of $C^{R}_{\ast}(X)$.
We thus have the quotient chain complex $C^{Q}_{\ast}(X) = C^{R}_{\ast}(X) / C^{D}_{\ast}(X)$.

Let $M$ be a right $\mathbb{Z}[G_{X}]$-module and $N$ a left $\mathbb{Z}[G_{X}]$-module.
The $n$-th \emph{quandle homology group} $H^{Q}_{\ast}(X; M)$ of $X$ with coefficient in $M$ is the $n$-th homology group of the chain complex $C^{Q}_{\ast}(X; M) = M \otimes_{\mathbb{Z}[G_{X}]} C^{Q}_{\ast}(X)$.
The $n$-th \emph{quandle cohomology group} $H_{Q}^{n}(X; N)$ of $X$ with coefficient in $N$ is the $n$-th cohomology group of the cochain complex $C_{Q}^{\ast}(X; N) = \mathrm{Hom}_{\mathbb{Z}[G_{X}]}(C^{Q}_{\ast}(X), N)$.
Here, $\mathrm{Hom}_{\mathbb{Z}[G_{X}]}(C^{Q}_{n}(X), N)$ denotes the abelian group consisting of $\mathbb{Z}[G_{X}]$-homomorphisms $C^{Q}_{n}(X) \rightarrow N$.
For the chain complex $C^{R}_{\ast}(X)$, we also define the \emph{rack homology group} $H^{R}_{\ast}(X; M)$ and the \emph{rack cohomology group} $H_{R}^{\ast}(X; N)$ in the same manner.

Let $A$ be an abelian group and $\mathrm{Hom}(M, A)$ denote the abelian group consisting of homomorphisms $M \rightarrow A$.
Then $\mathrm{Hom}(M, A)$ is a left $\mathbb{Z}[G_{X}]$-module by $\alpha f(r) = f(r \alpha)$ for any $f \in \mathrm{Hom}(M, A)$, $\alpha \in \mathbb{Z}[G_{X}]$, and $r \in M$.
We thus have the rack or quandle cohomology groups $H_{W}^{n}(X; \mathrm{Hom}(M, A))$, where the letter $W$ stands for $R$ or $Q$.
Since $C_{W}^{n}(X; \mathrm{Hom}(M, A))$ is isomorphic to $\mathrm{Hom}(C^{W}_{n}(X; M), A)$, we can define a pairing $\langle \ , \ \rangle : C_{W}^{n}(X; \mathrm{Hom}(M, A)) \otimes C^{W}_{n}(X; M) \rightarrow A$ by
\[
 \langle f, r \otimes (x_{1}, x_{2}, \cdots, x_{n}) \rangle = f(x_{1}, x_{2}, \cdots, x_{n})(r).
\]
We thus have a pairing $\langle \ , \ \rangle : H_{W}^{n}(X; \mathrm{Hom}(M, A)) \otimes H^{W}_{n}(X; M) \rightarrow A$.

Let $Y$ be a set equipped with a right action of $G_{X}$.
Then the free abelian group $\mathbb{Z}[Y]$ is a right $\mathbb{Z}[G_{X}]$-module.
We thus have the rack or quandle homology group $H^{W}_{n}(X; \mathbb{Z}[Y])$, the cohomology group $H_{W}^{n}(X; \mathrm{Hom}(\mathbb{Z}[Y], A))$, and a pairing $\langle \ , \ \rangle : H_{W}^{n}(X; \mathrm{Hom}(\mathbb{Z}[Y], A)) \otimes H^{W}_{n}(X; \mathbb{Z}[Y]) \rightarrow A$.

\begin{remark}
Our quandle homology group is isomorphic to the quandle homology group in \cite{CEGS99} if we regard a $\mathbb{Z}[G_{X}]$-module as a $\mathbb{Z}(X)$-module by $\eta_{x, y}(\alpha) = y^{-1} \alpha$ and $\tau_{x, y}(\alpha) = (1 - (x \ast y)^{-1}) \alpha$ for any $x, y \in X$ and $\alpha \in \mathbb{Z}[G_{X}]$.
A homomorphism from our $C^{Q}_{n}(X)$ to $C_{n}$ in \cite{CEGS99} sending $(x_{1}, x_{2}, \cdots, x_{n})$ to $\pm x_{1} x_{2} \cdots x_{n} (x_{1}, x_{2}, \cdots, x_{n})$ induces an isomorphism.
\end{remark}

\begin{remark}
Our quandle homology group $H^{Q}_{n}(X; \mathbb{Z}[Y])$ is isomorphic to a quandle homology group $H^{Q}_{n}(X)_{Y}$ in \cite{Kamada06}.
A homomorphism from our $C^{Q}_{n}(X; \mathbb{Z}[Y])$ to $C^{Q}_{n}(X)_{Y}$ in \cite{Kamada06} sending $r \otimes (x_{1}, x_{2}, \cdots, x_{n})$ to $(r, x_{1}, x_{2}, \cdots, x_{n})$ induces an isomorphism.
\end{remark}

\section{Simplicial quandle homology}
\label{sec_main}
In this section, we introduce simplicial quandle homology.
We show that we can construct a homomorphism from an $n$-th rack or quandle homology group to an $(n+1)$-th simplicial quandle homology group.

\subsection{Simplicial quandle homology}
Let $X$ be a quandle.
Consider the free abelian group $C^{\Delta}_{n}(X)$ generated by all $(n+1)$-tuples $(x_{0}, x_{1}, \cdots, x_{n}) \in X^{n+1}$ for each $n \geq 0$.
Define a map $\partial : C^{\Delta}_{n}(X) \rightarrow C^{\Delta}_{n-1}(X)$ by
\[
 \partial(x_{0}, x_{1}, \cdots, x_{n}) = \sum_{i = 0}^{n} (-1)^{i} (x_{0}, x_{1}, \cdots, \widehat{x_i}, \cdots, x_{n}).
\]
Then $\partial$ satisfies $\partial \circ \partial = 0$.
Thus, $C^{\Delta}_{\ast}(X) = (C^{\Delta}_{n}(X), \partial)$ is a chain complex.
It is easy to see that $C^{\Delta}_{\ast}(X)$ is acyclic.
A generator of $C^{\Delta}_{n}(X)$ may be identified with an $n$-simplex whose vertices are labelled by elements of $X$.
The boundary map $\partial$ sends a simplex to a formal sum of its boundaries.

The associated group $G_{X}$ acts on $C^{\Delta}_{n}(X)$ from the right by 
$(x_{0}, x_{1}, \cdots, x_{n}) \hskip 0.1em g = (x_{0} \ast g, x_{1} \ast g, \cdots, x_{n} \ast g)$ 
for each $(x_{0}, x_{1}, \cdots, x_{n}) \in C^{\Delta}_{n}(X)$ and $g \in G_{X}$.
Thus, $C^{\Delta}_{n}(X)$ is a right $\mathbb{Z}[G_{X}]$-module.
We let $C^{\Delta}_{n}(X)_{G_{X}} = C^{\Delta}_{n}(X) \otimes_{\mathbb{Z}[G_{X}]} \mathbb{Z}$.
Then $C^{\Delta}_{\ast}(X)_{G_{X}} = (C^{\Delta}_{n}(X)_{G_{X}}, \partial)$ is obviously a chain complex.

The $n$-th \emph{simplicial quandle homology group} $H^{\Delta}_{n}(X)$ of $X$ is the $n$-th homology group of the chain complex $C^{\Delta}_{\ast}(X)_{G_{X}}$.
Let $A$ be an abelian group.
The $n$-th \emph{simplicial quandle cohomology group} $H_{\Delta}^{n}(X; A)$ of $X$ with coefficient in $A$ is the $n$-th cohomology group of 
the cochain complex $\mathrm{Hom}_{\mathbb{Z}[G_{X}]}(C^{\Delta}_{\ast}(X), A)$.

\subsection{Quandle homology and simplicial quandle homology}
Recall that the associated group $G_{X}$ acts on $X$ from the right.
We thus have the rack or quandle homology group $H^{W}_{\ast}(X; \mathbb{Z}[X])$ ($W = R, Q$).

Let $I_{n}$ be the set consisting of maps $\iota : \{ 1, 2, \cdots, n\} \rightarrow \{ 0, 1 \}$.
Associated with $r \otimes (x_{1}, x_{2}, \cdots, x_{n}) \in C^{R}_{n}(X; \mathbb{Z}[X])$ and $\iota \in I_{n}$, define elements $r(\iota), x(\iota, i) \in X$ by
\begin{eqnarray*}
\begin{split}
 r(\iota) & = r \ast (x_{1}^{\iota(1)} x_{2}^{\iota(2)} \cdots x_{n}^{\iota(n)}), \\
 x(\iota,i) & = x_{i} \ast (x_{i+1}^{\iota(i+1)} x_{i+2}^{\iota(i+2)} \cdots x_{n}^{\iota(n)}).
\end{split}
\end{eqnarray*}
Suppose $|\iota|$ denotes the cardinality of the set $\{ i \mid \iota(i) = 1, \, 1 \leq i \leq n \}$.

\begin{figure}[htb]
\begin{center}
\includegraphics[scale=0.25]{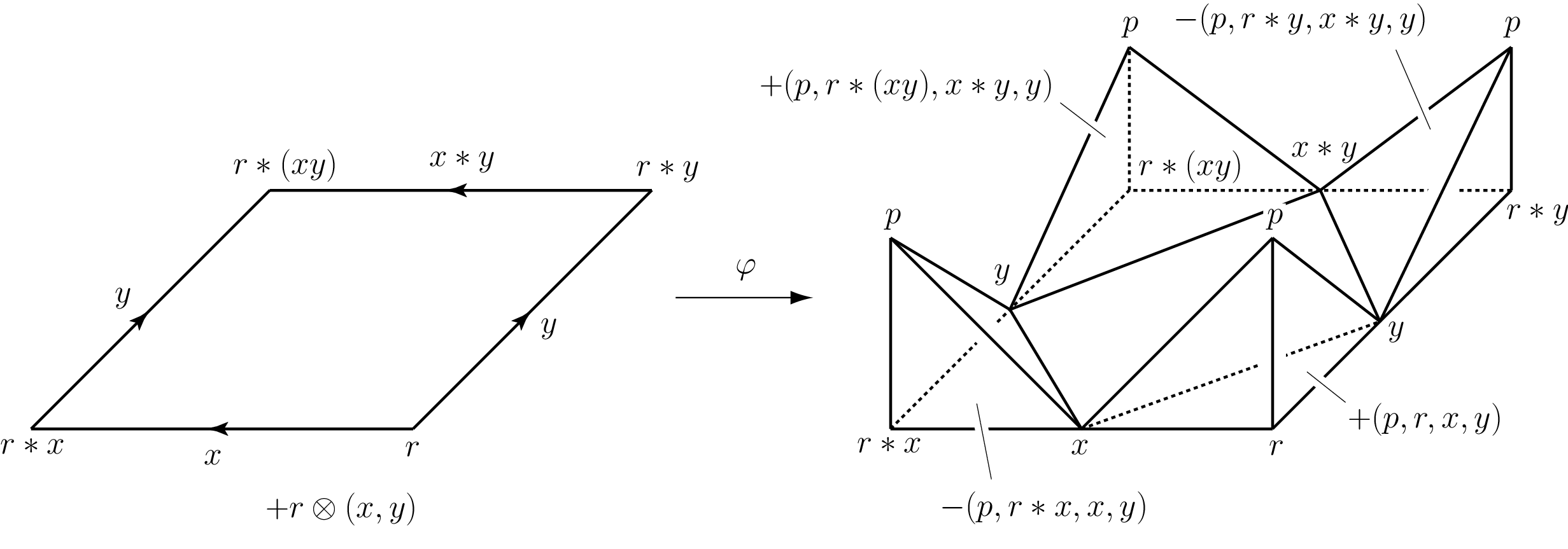}
\end{center}
\caption{A graphical explanation of the homomorphism $\varphi$.}
\label{fig:taking_cone_positive}
\end{figure}

Choose and fix an element $p \in X$.
For each $n \geq 0$, define a homomorphism $\varphi : C^{R}_{n}(X; \mathbb{Z}[X]) \rightarrow C^{\Delta}_{n+1}(X)_{G_{X}}$ by
\begin{equation*}
\varphi( r \otimes (x_{1}, x_{2}, \cdots, x_{n})) = \sum_{\iota \in I_{n}} (-1)^{|\iota|} (p, r(\iota), x(\iota,1), x(\iota,2), \cdots, x(\iota,n)).
\end{equation*}
In particular, in $n = 2$ case,
\[
\begin{split}
 \varphi(r \otimes (x, y)) = (p, r, x, y) - (p, r \ast x, x, y) - (p, r \ast y, x \ast y, y) + (p, r \ast (xy), x \ast y, y)
\end{split}
\]
(see also Figure \ref{fig:taking_cone_positive}).
Then we have the following lemma.

\begin{lemma}
\label{lem:chain_map}
The homomorphism $\varphi : C^{R}_{n}(X; \mathbb{Z}[X]) \rightarrow C^{\Delta}_{n+1}(X)_{G_{X}}$ is a chain map.
\end{lemma}

\begin{proof}
We let $I_{n}^{(i)} = \{ \iota \in I_{n} \mid \iota(i) = 0 \}$, $1 \leq i \leq n$.
For each $\iota \in I_{n}^{(i)}$, define a map $\iota^{(i)} \in I_{n} \setminus I_{n}^{(i)}$ by
\[
\iota^{(i)}(j) = 
 \begin{cases}
  1 \ & \mathrm{if} \ j = i, \\
  \iota(j) & \mathrm{if} \ j \neq i.
 \end{cases}
\]
Then $I_{n}$ is obviously decomposed into $I_{n}^{(i)}$ and a set $\{\iota^{(i)} \mid \iota \in I_{n}^{(i)} \}$.
For any $\iota \in I_{n}^{(i)}$, we have $|\iota^{(i)}| = |\iota| + 1$.

One computes
\begin{eqnarray*}
\begin{split}
 \partial & (\varphi( r \otimes (x_{1}, \cdots, x_{n}))) \\
 & = \sum_{\iota \in I_{n}} (-1)^{|\iota|} (\widehat{p}, r(\iota), x(\iota, 1), \cdots, x(\iota, n))
   - \sum_{\iota \in I_{n}} (-1)^{|\iota|} (p, \widehat{r(\iota)}, x(\iota, 1), \cdots, x(\iota, n)) \\
 & \quad + \sum_{\iota \in I_{n}} (-1)^{|\iota|} \sum_{i = 1}^{n} (-1)^{i} (p, r(\iota), x(\iota, 1), \cdots, \widehat{x(\iota, i)}, \cdots, x(\iota,n)) \\
 & = \sum_{\iota \in I_{n}^{(n)}} (-1)^{|\iota|} \{ (r(\iota), x(\iota, 1), \cdots, x(\iota, n))
   - (r(\iota^{(n)}), x(\iota^{(n)}, 1), \cdots, x(\iota^{(n)}, n)) \} \\
 & \quad - \sum_{\iota \in I_{n}^{(1)}} (-1)^{|\iota|} \{ (p, x(\iota, 1), \cdots, x(\iota, n))
         - (p, x(\iota^{(1)}, 1), \cdots, x(\iota^{(1)}, n)) \} \\
 & \qquad + \sum_{\iota \in I_{n}} (-1)^{|\iota|} \sum_{i = 1}^{n} (-1)^{i} (p, r(\iota), x(\iota, 1), \cdots, \widehat{x(\iota, i)}, \cdots, x(\iota, n)).
\end{split}
\end{eqnarray*}
By definition,
\begin{eqnarray*}
\begin{split}
 (r(\iota^{(n)}), x(\iota^{(n)}, 1), \cdots, x(\iota^{(n)}, n))
 & = (r(\iota) \ast x_{n}, x(\iota, 1) \ast x_{n}, \cdots, x(\iota, n) \ast x_{n}) \\
 & = (r(\iota), x(\iota, 1), \cdots, x(\iota, n)) \hskip 0.2em x_{n}
\end{split}
\end{eqnarray*}
for any $\iota \in I_{n}^{(n)}$.
Thus, we have
\begin{equation}
\label{eq:chain_map_1}
 (r(\iota^{(n)}), x(\iota^{(n)}, 1), \cdots, x(\iota^{(n)}, n)) = (r(\iota), x(\iota, 1), \cdots, x(\iota, n))
\end{equation}
in $C^{\Delta}_{n}(X)_{G_{X}}$.
Further, since $x(\iota, i)$ does not depend on the value of $\iota(1)$, we have
\begin{equation}
\label{eq:chain_map_2}
 x(\iota, i) = x(\iota^{(1)}, i)
\end{equation}
for any $\iota \in I_{n}^{(1)}$ and any index $i$.
By (\ref{eq:chain_map_1}) and (\ref{eq:chain_map_2}),
\begin{equation}
\label{eq:boundary_first}
\begin{split}
 & \partial (\varphi( r \otimes (x_{1}, \cdots, x_{n}))) = \\
 & \qquad \sum_{\iota \in I_{n}} (-1)^{|\iota|} \sum_{i = 1}^{n} (-1)^{i} (p, r(\iota), x(\iota, 1), \cdots, \widehat{x(\iota, i)}, \cdots, x(\iota, n)). \\
\end{split}
\end{equation}
On the other hand,
\begin{equation}
\label{eq:boundary_second}
\begin{split}
 \varphi(\partial & (r \otimes (x_{1}, \cdots, x_{n}))) \\
 & = \varphi \biggl( \sum_{i=1}^{n} (-1)^{i} \{ r \otimes (x_{1}, \cdots, \widehat{x_{i}}, \cdots, x_{n}) \\
 & \qquad - r \ast x_{i} \otimes (x_{1} \ast x_{i}, \cdots, x_{i-1} \ast x_{i}, x_{i+1}, \cdots, x_{n}) \} \biggr) \\
 & = \sum_{i=1}^{n} (-1)^{i} \sum_{\iota \in I_{n}^{(i)}} (-1)^{|\iota|}
     \{ (p, r(\iota), x(\iota, 1), \cdots, \widehat{x(\iota, i)}, \cdots, x(\iota, n)) \\
 & \qquad - (p, r(\iota^{(i)}), x(\iota^{(i)}, 1), \cdots, \widehat{x(\iota^{(i)}, i)}, \cdots, x(\iota^{(i)}, n)) \} \\
 & = \sum_{\iota \in I_{n}} (-1)^{|\iota|} \sum_{i = 1}^{n} (-1)^{i} (p, r(\iota), x(\iota, 1), \cdots, \widehat{x(\iota, i)}, \cdots, x(\iota, n)).
\end{split}
\end{equation}
By (\ref{eq:boundary_first}) and (\ref{eq:boundary_second}), we obtain $\partial ( \varphi(r \otimes (x_{1}, \cdots, x_{n}))) = \varphi( \partial(r \otimes (x_{1}, \cdots, x_{n})))$.
\end{proof}

Lemma \ref{lem:chain_map} immediately gives us the following theorem.

\begin{theorem}
\label{thm:homomorphisms}
For any quandle $X$ and each $n \geq 0$, there is a homomorphism
\[
 \varphi : H^{R}_{n}(X; \mathbb{Z}[X]) \longrightarrow H^{\Delta}_{n+1}(X).
\]
\end{theorem}

\begin{proposition}
\label{prop:does_not_depend}
The homomorphism $\varphi : H^{R}_{n}(X; \mathbb{Z}[X]) \rightarrow H^{\Delta}_{n+1}(X)$ does not depend on the choice of an element $p \in X$.
\end{proposition}

\begin{proof}
Let $p$ and $p^{\prime}$ be elements of $X$.
Suppose $\varphi$ and $\varphi^{\prime}$ are the chain maps related to $p$ and $p^{\prime}$ respectively.
Define a homomorphism $\Phi: C^{R}_{n}(X, \mathbb{Z}[X]) \rightarrow C^{\Delta}_{n+2}(X)_{G_{X}}$ by
\[
 \Phi(r \otimes (x_{1}, \cdots, x_{n})) = \sum_{\iota \in I_{n}} (-1)^{|\iota|} (p^{\prime}, p, r(\iota), x(\iota, 1), \cdots, x(\iota, n)).
\]
Then it is routine to check that $\Phi \circ \partial + \partial \circ \Phi = \varphi - \varphi^{\prime}$.
That is, $\Phi$ is a chain homotopy between $\varphi$ and $\varphi^{\prime}$.
\end{proof}

\begin{theorem}
\label{thm:homomorphisms_for_2}
For any quandle $X$, there is a homomorphism
\[
 \varphi : H^{Q}_{2}(X; \mathbb{Z}[X]) \longrightarrow H^{\Delta}_{3}(X).
\]
\end{theorem}
\begin{proof}
For any $(x_{0}, x_{1}, x_{2}, x_{3}, x_{4}) \in C^{\Delta}_{4}(X)_{G_{X}}$, we have
\begin{eqnarray*}
\begin{split}
 \qquad & \hskip -2em (x_{0}, x_{2}, x_{3}, x_{4}) - (x_{0}, x_{1}, x_{3}, x_{4}) \\
 \qquad = & \> + (x_{0}, x_{1}, x_{2}, x_{3}) - (x_{0}, x_{1}, x_{2}, x_{4}) + (x_{1}, x_{2}, x_{3}, x_{4}) - \partial(x_{0}, x_{1}, x_{2}, x_{3}, x_{4}).
\end{split}
\end{eqnarray*}
Hence,
\begin{eqnarray}
\label{eq:does_not_degenerate}
\begin{split}
 \qquad & \hskip -2em \varphi(r \otimes (x, y)) \\
 \qquad = & \> (p, r, x, y) - (p, r \ast x, x, y) - (p, r \ast y, x \ast y, y) + (p, r \ast (xy), x \ast y, y) \\
 \qquad = & \> (p, r \ast x, r, x) - (p, r \ast x, r, y) + (r \ast x, r, x, y) - \partial(p, r \ast x, r, x, y) \\
 \qquad & \> \enskip - (p, r \ast (xy), r \ast y, x \ast y) + (p, r \ast (xy), r \ast y, y) \\
 \qquad & \> \enskip \enskip - (r \ast (xy), r \ast y, x \ast y, y) + \partial(p, r \ast (xy), r \ast y, x \ast y, y) \\
 \qquad = & \> (p, r \ast x, r, x) - (p, r \ast x, r, y) \\
 \qquad & \> \enskip - (p, r \ast (xy), r \ast y, x \ast y) + (p, r \ast (xy), r \ast y, y) \\
 \qquad & \> \enskip \enskip - \partial(p, r \ast x, r, x, y) + \partial(p, r \ast (xy), r \ast y, x \ast y, y)
\end{split}
\end{eqnarray}
for each $r \otimes (x, y) \in C^{R}_{2}(X; \mathbb{Z}[X])$.
Here, the third equality follows from equations $(r \ast (xy), r \ast y, x \ast y, y) = (r \ast x, r, x, y) y = (r \ast x, r, x, y)$ in $C^{\Delta}_{3}(X)_{G_{X}}$.
Thus, for any $r \otimes (x, x) \in C^{D}_{2}(X; \mathbb{Z}[X])$, we have
\[
\varphi(r \otimes (x, x)) = - \partial(p, r \ast x, r, x, x) + \partial(p, r \ast x^{2}, r \ast x, x, x).
\]
\end{proof}

\section{Shadow coloring and fundamental class}
\label{sec_shadow_coloring}
In this section, we recall the definition of a shadow coloring.
It is known that, associated with a shadow coloring, we have a class in a second quandle homology group.
This homology class is called the fundamental class of a shadow coloring.
We show that a fundamental class derived from a shadow coloring is determined by the ``conjugacy class'' of the shadow coloring.

\subsection{Shadow coloring}
\label{subsec:def_of_shadow_colorings}
\begin{figure}[htb]
\begin{center}
\includegraphics[scale=0.32]{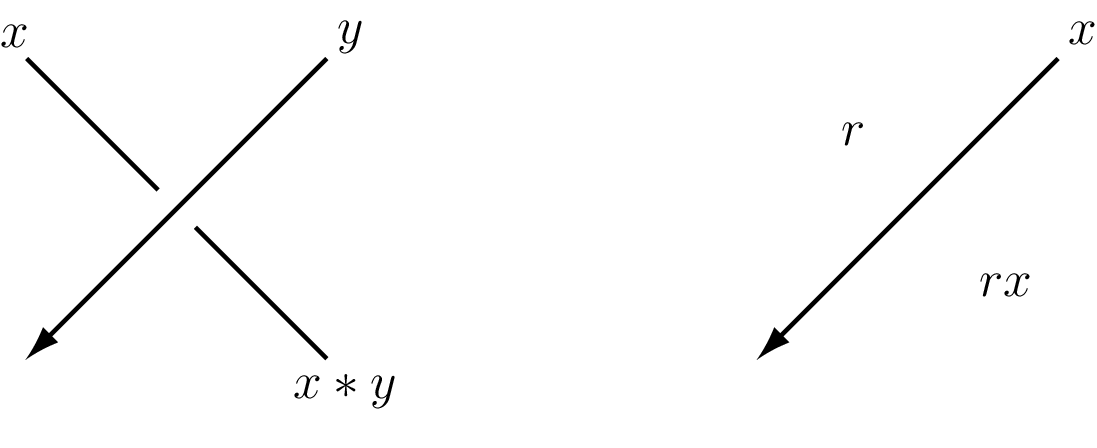}
\end{center}
\caption{Rules for colorings.}
\label{fig:coloring}
\end{figure}

Let $X$ be a quandle and $Y$ a set equipped with a right action of the associated group $G_{X}$.
Suppose $L$ is an oriented link in $S^{3}$ and $D$ a diagram of $L$.
An {\it arc coloring} of $D$ is a map $\mathcal{A} : \{ \textrm{arcs of} \ D \} \rightarrow X$ satisfying the condition 
illustrated in the left-hand side of Figure \ref{fig:coloring} at each crossing.
If an arc coloring of $D$ is assigned, a {\it region coloring} of $D$ is a map $\mathcal{R} : \{ \textrm{regions of} \ D \} \rightarrow Y$ satisfying the condition 
illustrated in the right-hand side of Figure \ref{fig:coloring} around each arc.
We call a pair $\mathcal{S} = (\mathcal{A}, \mathcal{R})$ a {\it shadow coloring} of $D$.

Let $m_{1}, m_{2}, \cdots, m_{n}$ be the Wirtinger generators of $\pi_{1}(S^{3} \setminus L)$ related to the arcs $\alpha_{1}, \alpha_{2}, \cdots, \alpha_{n}$ of $D$ respectively.
Then, associated with an arc coloring $\mathcal{A}$, we obtain a representation $\rho_{\mathcal{A}} : \pi_{1}(S^{3} \setminus L) \rightarrow G_{X}$ 
which sends each $m_{i}$ to $\mathcal{A}(\alpha_{i})$.

\subsection{Fundamental class}
\begin{figure}[htb]
\begin{center}
\includegraphics[scale=0.30]{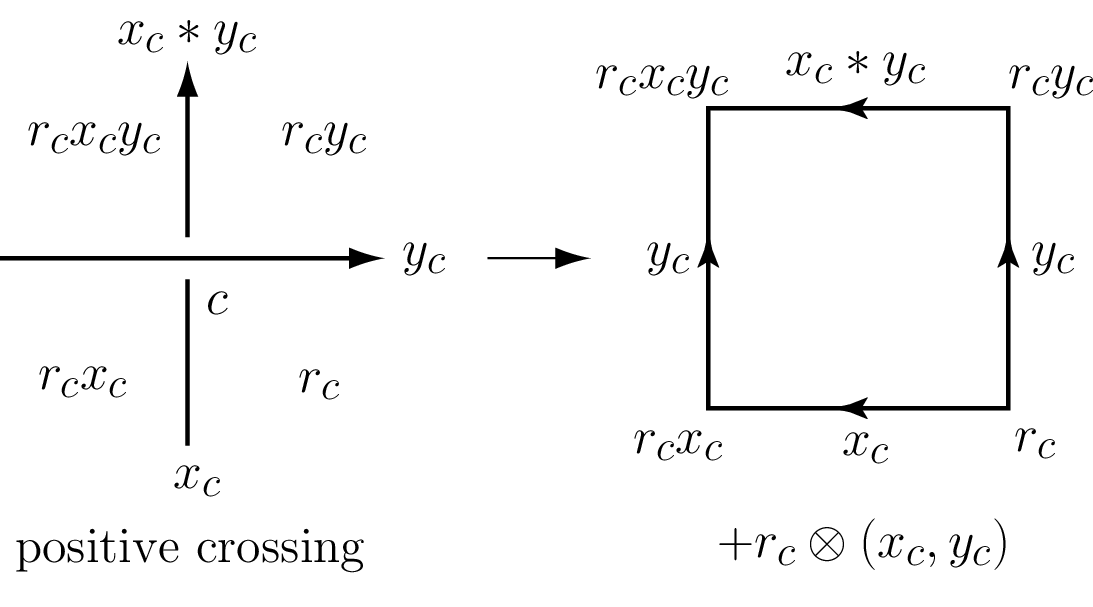}
\qquad
\includegraphics[scale=0.30]{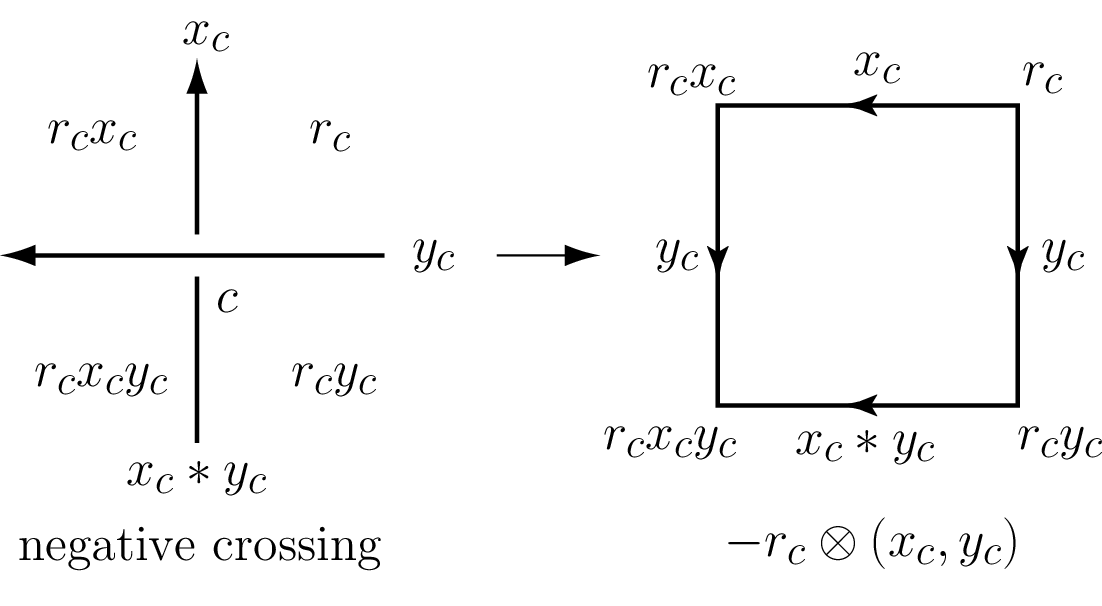}
\end{center}
\caption{We assign a square for each crossing.}
\label{fig:cycle_from_shadow_coloring}
\end{figure}

For a shadow coloring $\mathcal{S}$, define a chain
\[
 C(\mathcal{S}) = \sum_{c} \varepsilon_{c} r_{c} \otimes (x_{c}, y_{c}) \in C^{Q}_{2}(X; \mathbb{Z}[Y]).
\]
Here, the sum runs over all crossings $c$ of $D$, $\varepsilon_{c}$ is $1$ or $-1$ depending on whether $c$ is positive or negative respectively, 
and $x_{c}, y_{c} \in X$ and $r_{c} \in Y$ denote colors around $c$ as depicted in Figure \ref{fig:cycle_from_shadow_coloring}.
It is straightforward to check that we have the following lemma.

\begin{lemma}[\cite{CKS01}, \cite{Kamada06}]
The chain $C(\mathcal{S})$ is a cycle.
\end{lemma}

Let $D^{\prime}$ be another diagram of $L$ obtained from $D$ by a single Reidemeister move.
Then there is a unique shadow coloring $\mathcal{S^{\prime}}$ of $D^{\prime}$ which coincides with $\mathcal{S}$ except for colors of arcs and regions related to the move.

\begin{lemma}[\cite{CKS01}, \cite{Kamada06}]
\label{lem:independent_to_diagrams}
The cycles $C(\mathcal{S})$ and $C(\mathcal{S}^{\prime})$ are homologous.
\end{lemma}

\begin{proof}
A first Reidemeister move adds or subtracts $\pm r \otimes (x, x)$ to or from $C(\mathcal{S})$ for some $r \in Y$ and $x \in X$, 
but $\pm r \otimes (x, x) = 0$ in $C^{Q}_{2}(X; \mathbb{Z}[Y])$. 
A second Reidemeister move adds or subtracts $r \otimes (x, y) - r \otimes (x, y) = 0$ for some $r \in Y$ and $x, y \in X$.
A third Reidemeister move adds $\pm \partial(r \otimes (x, y, z))$ to $C(\mathcal{S})$ for some $r \in Y$ and $x, y, z \in X$ (see Figure \ref{fig:reidemeister_III}), 
therefore it does not change the homology class.
\begin{figure}[htb]
\begin{center}
\includegraphics[scale=0.32]{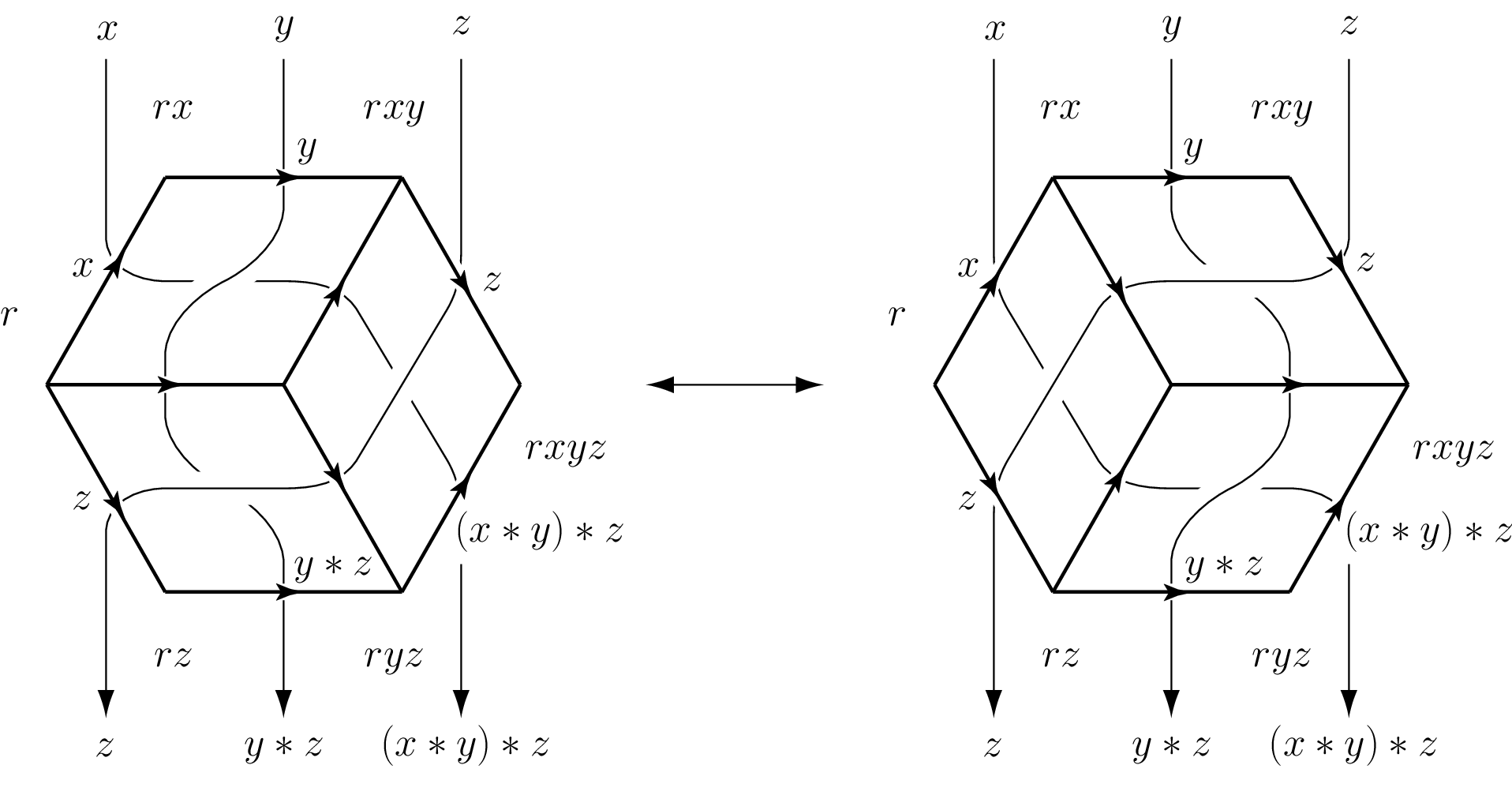}
\end{center}
\caption{The effect of a third Reidemeister move is identified with an addition of the boundaries of a $3$-cube.}
\label{fig:reidemeister_III}
\end{figure}
\end{proof}

We call the homology class $[C(\mathcal{S})] \in H^{Q}_{2}(X; \mathbb{Z}[Y])$ a \emph{fundamental class} derived from $\mathcal{S}$.

\subsection{}
We devote the remaining of this section to show the following theorem.

\begin{theorem}
\label{thm:region_and_conjugacy_invariance}
Let $X$ be a quandle, $G_{X}$ the associated group of $X$, and $Y$ a set equipped with a right action of $G_{X}$.
Suppose $L$ is an oriented link in $S^{3}$, $D$ a diagram of $L$, and $\mathcal{S} = (\mathcal{A}, \mathcal{R})$ a shadow coloring of $D$ with respect to $X$ and $Y$.
If $G_{X}$ acts on $Y$ transitively, then a fundamental class $[C(\mathcal{S})]$ derived from $\mathcal{S}$ does not depend on the choice of $\mathcal{R}$.
Further, if the natural map $X \to G_X$ is injective, then $[C(\mathcal{S})]$ is determined by the conjugacy class of 
the representation $\rho_{\mathcal{A}} : \pi_{1}(S^{3} \setminus L) \rightarrow G_{X}$ derived from $\mathcal{A}$.
\end{theorem}

To prove Theorem \ref{thm:region_and_conjugacy_invariance}, we first show the following lemma.
For any $r \in Y$, let $\mathcal{R}_{r}$ denote a region coloring of $D$ sending the region containing the point at infinity to $r$.
Remark that $\mathcal{R}_{r}$ surely exists and is unique.

\begin{lemma}
\label{lem:independent_to_region_colorings}
For any $r \in Y$ and $w \in X$, $C((\mathcal{A}, \mathcal{R}_{r}))$ and $C((\mathcal{A}, \mathcal{R}_{r w}))$ are homologous.
\end{lemma}

\begin{proof}
We first consider a region coloring of $D$ with respect to $G_{X}$ instead of $Y$.
For any $h \in G_{X}$, let $\mathcal{R}^{G_{X}}_{h}$ denote a region coloring of $D$ with respect to $G_{X}$ sending the region containing the point at infinity to $h$.
Assume that $C((\mathcal{A}, \mathcal{R}^{G_{X}}_{1})) = \sum_{c} \varepsilon_{c} g_{c} \otimes (x_{c}, y_{c})$ for some $g_{c} \in G_{X}$.
Here, $x_{c}, y_{c} \in X$ denote colors around a crossing $c$ of $D$ with respect to $\mathcal{A}$.
Then it is routine to check that $C((\mathcal{A}, \mathcal{R}^{G_{X}}_{h})) = \sum_{c} \varepsilon_{c} h g_{c} \otimes (x_{c}, y_{c})$.

Define a homomorphism $\zeta : C^{Q}_{1}(X, \mathbb{Z}[G_{X}]) \rightarrow C^{Q}_{2}(X, \mathbb{Z}[G_{X}])$ by $\zeta(g \otimes (x) ) = g \otimes (w \ast g, x)$.
Then
\begin{equation}
\label{eq:universal_case}
\begin{split}
   & \partial \left( \sum_{c} \varepsilon_{c} g_{c} \otimes (w \ast g_{c}, x_{c}, y_{c}) \right) \\
 = & \sum_{c} \varepsilon_{c} \{ - \, g_{c} \otimes (x_{c}, y_{c}) +  g_{c} (w \ast g_{c}) \otimes (x_{c}, y_{c}) \\
   & \hskip 2.9em \quad + g_{c} \otimes (w \ast g_{c}, y_{c}) - g_{c} x_{c} \otimes ((w \ast g_{c}) \ast x_{c}, y_{c}) \\
   & \hskip 2.9em \quad \quad - g_{c} \otimes (w \ast g_{c}, x_{c}) + g_{c} y_{c} \otimes ((w \ast g_{c}) \ast y_{c}, x_{c} \ast y_{c}) \} \\
 = & - C((\mathcal{A}, \mathcal{R}^{G_{X}}_{1})) + C((\mathcal{A}, \mathcal{R}^{G_{X}}_{w})) \\
   & \quad + \sum_c \varepsilon_{c} \{ \zeta ( g_{c} \otimes (y_{c}) - g_{c} x_{c} \otimes (y_{c}) - g_{c} \otimes (x_{c}) + g_{c} y_{c} \otimes (x_{c} \ast y_{c})) \} \\
 = & - C((\mathcal{A}, \mathcal{R}^{G_{X}}_{1})) + C((\mathcal{A}, \mathcal{R}^{G_{X}}_{w})) - \zeta ( \partial ( \sum_{c} \varepsilon_{c} g_{c} \otimes (x_{c}, y_{c})) ) \\
 = & - C((\mathcal{A}, \mathcal{R}^{G_{X}}_{1})) + C((\mathcal{A}, \mathcal{R}^{G_{X}}_{w})) - \zeta ( \partial C((\mathcal{A}, \mathcal{R}^{G_{X}}_{1})) ) \\
 = & - C((\mathcal{A}, \mathcal{R}^{G_{X}}_{1})) + C((\mathcal{A}, \mathcal{R}^{G_{X}}_{w})).
\end{split}
\end{equation}
Here, the second equality follows from the relations $g_{c} (w \ast g_{c}) = w g_{c}$, $(w \ast g_{c}) \ast x_{c} = w \ast (g_{c} x_{c})$, 
and $(w \ast g_{c}) \ast y_{c} = w \ast (g_{c} y_{c})$.
The last equality follows since $C((\mathcal{A}, \mathcal{R}^{G_{X}}_{1}))$ is a cycle.

Define a chain map $\eta : C^{Q}_{n}(X, \mathbb{Z}[G_{X}]) \rightarrow C^{Q}_{n}(X, \mathbb{Z}[Y])$ by 
$\eta (g \otimes (x_{1}, x_{2}, \cdots, x_{n})) = r g \otimes (x_{1}, x_{2}, \cdots, x_{n})$.
Remark that $\eta (C((\mathcal{A}, \mathcal{R}^{G_{X}}_{h}))) = C((\mathcal{A}, \mathcal{R}_{r h}))$.
Applying $\eta$ to (\ref{eq:universal_case}), we obtain
\[
 \partial \left( \sum_{c} \varepsilon_{c} r g_{c} \otimes (w \ast g_{c}, x_{c}, y_{c}) \right) = - \> C((\mathcal{A}, \mathcal{R}_{r})) + C((\mathcal{A}, \mathcal{R}_{r w})).
\]
Therefore, $C((\mathcal{A}, \mathcal{R}_{r}))$ and $C((\mathcal{A}, \mathcal{R}_{r w}))$ are homologous.
\end{proof}

Let $\alpha_{1}, \alpha_{2}, \cdots, \alpha_{n}$ again be the arcs of $D$, and $\beta_{1}, \beta_{2}, \cdots, \beta_{m}$ be the regions of $D$.
It is easy to see that, for any arc coloring $\mathcal{A}$ and $g \in G_{X}$, 
a map $\mathcal{A} g : \{ \textrm{arcs of} \ D \} \rightarrow X$ sending $\alpha_{i}$ to $\mathcal{A}(\alpha_{i}) \ast g$ is also an arc coloring.
Similarly, for any region coloring $\mathcal{R}$, a map $\mathcal{R} g : \{ \textrm{regions of} \ D \} \rightarrow Y$ sending $\beta_{j}$ 
to $\mathcal{R}(\beta_{j}) g$ is also a region coloring.

\begin{lemma}
\label{lem:independent_to_conjugations}
For any $w \in X$, $C((\mathcal{A}, \mathcal{R}))$ and $C((\mathcal{A} w, \mathcal{R} w))$ are homologous.
\end{lemma}

\begin{proof}
Define a homomorphism $\theta : C^{Q}_{1}(X, \mathbb{Z}[Y]) \to C^{Q}_{2}(X, \mathbb{Z}[Y])$ by $\theta (r \otimes (x)) = r \otimes (x, w)$.
Suppose $x_{c}, y_{c} \in X$ and $r_{c} \in Y$ are colors around a crossing $c$ of $D$ with respect to a shadow coloring $(\mathcal{A}, \mathcal{R})$.
Then
\begin{eqnarray*}
\begin{split}
   & \partial \left( \sum_{c} \varepsilon_{c} r_{c} \otimes (x_{c}, y_{c}, w) \right) \\
 = & \sum_{c} \varepsilon_{c} \{ - r_{c} \otimes (y_{c}, w) +  r_{c} x_{c} \otimes (y_{c}, w) + r_{c} \otimes (x_{c}, w) - r_{c} y_{c} \otimes (x_{c} \ast y_{c}, w) \\
   & \quad - r_{c} \otimes (x_{c}, y_{c}) + r_{c} w \otimes (x_{c} \ast w, y_{c} \ast w) \} \\
 = & \> \theta \left( \sum_{c} \varepsilon_{c} \{ - r_{c} \otimes (y_{c}) + r_{c} x_{c} \otimes (y_{c}) + r_{c} \otimes (x_{c}) - r_{c} y_{c} \otimes (x_{c} \ast y_{c}) \} \right) \\
   & \quad - C((\mathcal{A}, \mathcal{R})) + C((\mathcal{A} w, \mathcal{R} w)) \\
 = & \> \theta ( \partial C((\mathcal{A}, \mathcal{R})) ) - C((\mathcal{A}, \mathcal{R})) + C((\mathcal{A} w, \mathcal{R} w)) \\
 = & - C((\mathcal{A}, \mathcal{R})) + C((\mathcal{A} w, \mathcal{R} w)).
\end{split}
\end{eqnarray*}
Here, the last equality follows since $C((\mathcal{A}, \mathcal{R}))$ is a cycle.
Therefore, $C((\mathcal{A}, \mathcal{R}))$ and $C((\mathcal{A} w, \mathcal{R} w))$ are homologous.
\end{proof}

\begin{proof}[Proof of Theorem \ref{thm:region_and_conjugacy_invariance}]
If $G_{X}$ acts on $Y$ transitively, then the homology class $[C(\mathcal{S})]$ does not depend on the choice of $\mathcal{R}$ by Lemma \ref{lem:independent_to_region_colorings}.

Let $\mathcal{A}$ and $\mathcal{A}'$ be two arc colorings such that $\rho_{\mathcal{A}'} = g^{-1} \rho_{\mathcal{A}} g$ for some $g \in G_X$.
Since $ g^{-1} \rho_{\mathcal{A}} g = \rho_{\mathcal{A} g}$, we have $\mathcal{A}' = \mathcal{A}g$ if the natural map $X \to G_X$ is injective.
On the other hand, $C((\mathcal{A}, \mathcal{R}))$ and $C((\mathcal{A} g, \mathcal{R} g))$ are homologous by Lemma \ref{lem:independent_to_conjugations}.
\end{proof}

\begin{remark}
\label{rem:general_cycles}
We have obtained a cycle $C(\mathcal{S}) \in C^{Q}_{2}(X, \mathbb{Z}[Y])$ derived from a shadow coloring $\mathcal{S}$.
Conversely, extending the work of Carter et al. \cite{CKS01}, we can show that, for any cycle $C \in C^{Q}_{2}(X, \mathbb{Z}[Y])$, 
there is a link diagram $D$ on an orientable surface and a shadow coloring $\mathcal{S}$ of $D$ with respect to $X$ and $Y$ such that $C(\mathcal{S}) = C$.
Further, we can show that a fundamental class $[C(\mathcal{S})] \in H^{Q}_{2}(X, \mathbb{Z}[Y])$ does not depend on the choice of a region coloring in the same line.
\end{remark}

\begin{remark}
Similar claims of Lemma \ref{lem:independent_to_conjugations} are proved by Etingof and Gra\~na \cite{etingof-grana}, 
and Niebrzydowski and Przytycki \cite{niebrzydowski-przytycki}.
\end{remark}

\section{Quandle consisting of all parabolic elements of $\PSLC$}
\label{sec_psl2c}
In this section, we define a quandle $\mathcal{P}$ consisting of all parabolic elements of $\PSLC$.
We see that $\mathcal{P}$ is identified with $\Cpm$.

\subsection{}
Let $\mathcal{P}$ be the set of all parabolic elements of $\PSLC$, then $\mathcal{P}$ is closed under conjugations.
Thus, $\mathcal{P}$ is a conjugation quandle with $x \ast y = y^{-1} x y$ for any $x, y \in \mathcal{P}$.
Let $G_{\mathcal{P}}$ be the associated group of $\mathcal{P}$.
The natural inclusion $\mathcal{P} \rightarrow \PSLC$ induces a homomorphism $\xi : G_{\mathcal{P}} \rightarrow \PSLC$.
Since $\PSLC$ is generated by parabolic elements, $\xi$ is surjective.

Let $L$ be an oriented link in $S^{3}$, $D$ a diagram of $L$, and $\mathcal{A}$ an arc coloring of $D$ with respect to $\mathcal{P}$.
Recall that we have a representation $\rho_{\mathcal{A}} : \pi_{1}(S^{3} \setminus L) \rightarrow G_{\mathcal{P}}$ derived from $\mathcal{A}$.
The composition $\xi \circ \rho_{\mathcal{A}} : \pi_{1}(S^{3} \setminus L) \rightarrow \PSLC$ sends meridians to parabolic elements.
We call a representation of $\pi_1(S^3 \setminus L)$ \emph{parabolic} if it sends each meridian to a parabolic element of $\PSLC$.
Thus $\xi \circ \rho_{\mathcal{A}}$ is a parabolic representation.
Conversely, let $\rho : \pi_{1}(S^{3} \setminus L) \rightarrow \PSLC$ be a parabolic representation.
We can define an arc coloring $\mathcal{A}_{\rho}$ satisfying $\xi \circ \rho_{\mathcal{A}_{\rho}} = \rho$.

The associated group $G_{\mathcal{P}}$ obviously acts on $\mathcal{P}$ transitively.
Further the first author \cite{Inoue10_2} showed that for any conjugation quandle $X$, the natural map $X \to G_X$ is injective.
Therefore  $[C(\mathcal{S})]$ is completely determined by the conjugacy class of $\xi \circ \rho_{\mathcal{A}}$.


\subsection{}
Let $\Cpm$ be the quotient of $\mathbb{C}^{2} \setminus \{ 0 \}$ by the equivalence relation $v \sim -v$.
Since each element of $\mathcal{P}$ has a presentation
\[
 \begin{pmatrix} a & b \\ c & d \end{pmatrix}^{-1} \begin{pmatrix} 1 & 1 \\ 0 & 1 \end{pmatrix} \begin{pmatrix} a & b \\ c & d \end{pmatrix} 
= \begin{pmatrix} 1 + c d & d^{2} \\ - \, c^{2} & 1 - c d \end{pmatrix}
\]
for some $\begin{pmatrix} a & b \\ c & d \end{pmatrix} \in \PSLC$, we can identify $\mathcal{P}$ with $\Cpm$ by 
\[
\begin{pmatrix} \alpha & \beta  \end{pmatrix} \leftrightarrow \begin{pmatrix} 1 + \alpha \beta & \beta^{2} \\ - \, \alpha^{2} & 1 - \alpha \beta \end{pmatrix}.
\]
The binary operation on $\mathcal{P}$ induces a binary operation on $\Cpm$ given by
\[
\begin{pmatrix} \alpha & \beta \end{pmatrix}  \ast \begin{pmatrix} \gamma & \delta \end{pmatrix} 
= \begin{pmatrix} \alpha & \beta \end{pmatrix} \begin{pmatrix} 1 + \gamma \delta & \delta^{2} \\ - \, \gamma^{2} & 1 - \gamma\delta \end{pmatrix}.
\]
Since $\PSLC$ usually acts on $\Cpm$ from left, we write it in the following form 
\begin{equation}
\label{eq:binary_operation_on_C^2}
\begin{pmatrix} \alpha \\ \beta \end{pmatrix}  \ast \begin{pmatrix} \gamma \\ \delta \end{pmatrix} 
= \begin{pmatrix} 1 + \gamma \delta & - \, \gamma^{2} \\ \delta^{2} & 1 - \gamma\delta \end{pmatrix} \begin{pmatrix} \alpha \\ \beta \end{pmatrix}.
\end{equation}
The inverse operation is given by 
\begin{equation}
\label{eq:inverse_binary_operation_on_C^2}
\begin{pmatrix} \alpha \\ \beta \end{pmatrix}  \ast^{-1} \begin{pmatrix} \gamma \\ \delta \end{pmatrix} 
= \begin{pmatrix} 1 - \gamma \delta & \gamma^{2} \\ - \, \delta^{2} & 1 + \gamma\delta \end{pmatrix} \begin{pmatrix} \alpha \\ \beta \end{pmatrix}.
\end{equation}

\begin{remark}
The set $\mathbb{C}^{2} \setminus \{ 0 \}$ is also a quandle with a binary operation given by (\ref{eq:binary_operation_on_C^2})
for $\begin{pmatrix} \alpha \\ \beta \end{pmatrix}, \begin{pmatrix} \gamma \\ \delta \end{pmatrix} \in \mathbb{C}^{2} \setminus \{ 0 \}$.
A natural projection $\mathbb{C}^{2} \setminus \{ 0 \} \rightarrow (\mathbb{C}^{2} \setminus \{ 0 \}) / \pm$ is a two-to-one quandle homomorphism.
\end{remark}

\section{Extended Bloch group}
\label{sec_extended_bloch}
In this section, we recall the definition of the \emph{extended Bloch group} $\hatBC$ by Neumann \cite{neumann04}.

\subsection{Bloch group}
Let $\mathbb{H}^3$ be the 3-dimensional hyperbolic space.
The Riemann sphere $\mathbb{C}P^1$ can be regarded as the ideal boundary of $\mathbb{H}^3$.
$\PSLC$ acts on $\mathbb{C}P^1$ by linear fractional transformations and it extends to an action on $\mathbb{H}^3$ as the group of orientation preserving isometries.
An \emph{ideal tetrahedron} is the convex hull of four distinct ordered points of $\mathbb{C}P^1$ in $\mathbb{H}^3$.
An ideal tetrahedron with ordered vertices  $z_0,z_1,z_2,z_3$ is parametrized by the cross ratio
\[
[z_0:z_1:z_2:z_3]=\frac{z_3-z_0}{z_3-z_1}\frac{z_2-z_1}{z_2-z_0}.
\]
The cross ratio satisfies $[g z_0:gz_1:gz_2:gz_3]=[z_0:z_1:z_2:z_3] $ for any $g \in \PSLC$. 
Let $[z_iz_j]$ be the edge of the ideal tetrahedron spanned by $z_i$ and $z_j$.
Take $\{i,j,k,l\}$ to be an even permutation of $\{0,1,2,3\}$.
We define the complex parameter $[z_{i}, z_{j}]$ of the edge by the cross ratio $[z_i:z_j:z_k:z_l]$.
This parameter only depends on the choice of the edge $[z_iz_j]$.
We can easily observe that the opposite edge has the same complex parameter.
If the complex parameter of the edge $[z_0z_1]$ (or $[z_2z_3]$) is $z$, 
then the complex parameter of the edge $[z_1z_2]$ (or $[z_0z_3]$) is $\frac{1}{1-z}$, 
and the complex parameter of the edge $[z_1z_3]$ (or $[z_0z_2]$) is $1-\frac{1}{z}$.

The pre-Bloch group $\mathcal{P}(\mathbb{C})$ is an abelian group generated by symbols $[z]$, $z \in \mathbb{C} \setminus \{0, 1 \}$, subject to the relation 
\[
[x]-[y]+[y/x]-\left[\frac{1-x^{-1}}{1-y^{-1}}\right]+\left[\frac{1-x}{1-y}\right] = 0, \quad x, y \in \mathbb{C} \setminus \{0,1\}.
\]
This relation is called the \emph{five term relation}.
This is equivalent to the relation $\sum_{i=0}^4 (-1)^i[z_0: \dots : \widehat{z_i} : \dots :z_4 ] =0$, where $z_0, \dots ,z_4$ are distinct five points of $\mathbb{C}P^1$.
We define a map $\lambda : \mathcal{P}(\mathbb{C}) \to \mathbb{C}^* \wedge_{\mathbb{Z}} \mathbb{C^*}$
by $[z] \to z \wedge (1-z)$.
The \emph{Bloch group} $\mathcal{B}(\mathbb{C})$ is the kernel of $\lambda$. 

\subsection{Extended Bloch group}
In this paper, we define the logarithm $\Log(z)$ by $\log|z| + i\arg(z)$ with $-\pi < \arg(z) \leq \pi$. 
Let $P$ be $\mathbb{C} \backslash \{ 0, 1\}$ cut along the lines $(-\infty, 0)$ and $(1, \infty)$.
There are two copies of these cut lines in $P$.
We denote a point on these lines by $x+0i$ (respectively $x-0i$) if the point on the boundary of upper (respectively lower)
half space.
We construct $\widehat{\mathbb{C}}$ by gluing $P \times \mathbb{Z} \times \mathbb{Z}$ along their boundaries by
\[
\begin{split}
(x+0i, p, q) &\sim (x-0i, p+2, q) \quad \textrm{for each $x \in (-\infty, 0)$}, \\
(x+0i, p, q) &\sim (x-0i, p, q+2) \quad \textrm{for each $x \in (1, \infty)$}.
\end{split}
\]
$\widehat{\mathbb{C}}$ consists of four components $X_{00}$, $X_{01}$, $X_{10}$ and $X_{11}$, 
where $X_{\varepsilon_1\varepsilon_2}$ is the component with $p \equiv \varepsilon_1$ and $q \equiv \varepsilon_2$ mod 2.
Each component is the universal abelian cover of $\mathbb{C} \backslash \{0,1\}$ and is considered as the Riemann surface of the multivalued function $(\log(z), -\log(1-z))$.
In fact, the map $(z; p,q) \mapsto (\Log(z)+ 2p\pi i, -\Log(1-z) + 2q \pi i)$ is a well-defined map from $\widehat{\mathbb{C}}$ to $\mathbb{C}^2$.

The extended pre-Bloch group is an abelian group generated by elements of $\widehat{\mathbb{C}}$ subject to some relations. 
To describe the relations, we need some definitions.

\begin{definition}
\label{def:comninatorial_flattenings}
Let $\Delta$ be an ideal tetrahedron with cross ratio $z$.
A \emph{combinatorial flattening} of $\Delta$ is a triple of complex numbers $(w_0,w_1,w_2)$ of the form 
\[
(w_0,w_1,w_2) = (\Log(z)+p\pi i, -\Log(1-z) + q\pi i, -\Log(z)+\Log(1-z) -p\pi i- q\pi i)
\]
for some $p, q \in \mathbb{Z}$.
\end{definition}
We call $w_0$ the \emph{log-parameter} of the edge $[z_0z_1]$ (or $[z_2z_3]$), 
$w_1$ of the edge $[z_1z_2]$ (or $[z_0z_3]$) and  
$w_2$ of the edge $[z_1z_3]$ (or $[z_0z_2]$).
We denote the log-parameter of an edge $E$ by $l_E$. 
We relate a combinatorial flattening with an element of $\widehat{\mathbb{C}}$ by the following map:
\[
(z;p,q) \mapsto (\Log(z)+ p \pi i, - \Log(1-z) + q \pi i, -\Log(z) + \Log(1-z) - p \pi i -q \pi i ).
\]
This map is a bijection (\cite{neumann04}, Lemma 3.2).

Consider the boundary of the ideal 4-simplex spanned by $z_0,z_1,z_2,z_3,z_4$.
Let $(w_0^i,w_1^i,w_2^i)$ be a combinatorial flattening of the ideal tetrahedron $(z_0, \dots , \widehat{z_i}, \dots, z_4)$.
The flattenings of these tetrahedra are said to satisfy \emph{flattening condition} if for each edge the signed sum of log parameters is zero.
There are ten edges $[z_iz_j]$ and each edge belongs to exactly three ideal tetrahedra.
The flattening condition is equivalent to the following ten equations:
\begin{equation}
\label{eq:lifted_five_term}
\begin{tabular}[]{lrclr}
$[z_0z_1]:$ & $ w_0^2-w_0^3+w_0^4=0$ & & $[z_0z_2]:$ & $-w_0^1-w_2^3+w_2^4=0$ \\
$[z_1z_2]:$ & $ w_0^0-w_1^3+w_1^4=0$ & & $[z_1z_3]:$ & $ w_2^0+w_1^2+w_2^4=0$ \\
$[z_2z_3]:$ & $ w_1^0-w_1^1+w_0^4=0$ & & $[z_2z_4]:$ & $ w_2^0-w_2^1-w_0^3=0$ \\
$[z_3z_4]:$ & $ w_0^0-w_0^1+w_0^2=0$ & & $[z_3z_0]:$ & $-w_2^1+w_2^2+w_1^4=0$ \\
$[z_4z_0]:$ & $-w_1^1+w_1^2-w_1^3=0$ & & $[z_4z_1]:$ & $ w_1^0-w_2^2-w_2^3=0$ 
\end{tabular}
\end{equation}

\begin{definition}
The \emph{extended pre-Bloch group} $\hatPC$ is an abelian group generated by elements of $\widehat{\mathbb{C}}$ subject to the following relations:
\begin{itemize}
\item[(i)]
$\sum_{i=0}^4 (-1)^i (w_0^i, w_1^i,w_2^i) = 0$, 
if the flattenings $(w_0^i,w_1^i,w_2^i)$ satisfy flattening condition.
\item[(ii)]
$
(z;p,q) + (z;p',q') = (z;p,q') + (z;p',q) \quad \textrm{with $p,p',q,q' \in \mathbb{Z}$.}
$
\end{itemize}
The first relation is called the \emph{lifted five term relation}.
The second relation is called the \emph{transfer relation}.
We denote the class of $(z; p, q)$ in $\hatPC$ by $[z; p, q]$.
\end{definition}

We define a map $\nu: \hatPC \to \mathbb{C} \wedge_{\mathbb{Z}} \mathbb{C} $ by
\[
[z;p,q] \mapsto (\Log(z)+p \pi i) \wedge (-\Log(1-z)+q \pi i).
\]
The kernel of $\nu$ is called the \emph{extended Bloch group} and denoted by $\hatBC$.

Let 
\[
R(z;p,q)=\mathcal{R}(z)+\frac{\pi i}{2} \left( q\Log(z)-p\Log \left(\frac{1}{1-z}\right) \right) - \frac{\pi^2}{6},
\]
where $\mathcal{R}(z)$ is given by
\[
\mathcal{R}(z) = - \int_{0}^{z} \frac{\Log(1-t)}{t} dt + \frac{1}{2}\Log(z)\Log(1-z).
\]
This map is well-defined on $\hatPC$ up to integer multiple of $\pi^2$.
So it induces a map $R: \hatPC \to \mathbb{C}/\pi^2\mathbb{Z}$.

\begin{theorem}[Neumann \cite{neumann04}, Theorem 12.1]
There exists an isomorphism $\lambda: H_3(\PSLC;\mathbb{Z}) \to \hatBC$ such that the composition $R \circ \lambda : H_3(\PSLC;\mathbb{Z}) \to \mathbb{C}/\pi^2\mathbb{Z}$
is the Cheeger-Chern-Simons class.
\end{theorem}

\section{Complex volume in terms of quandle homology}
\label{sec_invariant}
Let $G=\PSLC$ and $\mathcal{P}$ be the quandle consisting of all parabolic elements of $\PSLC$.
In this section, we construct a homomorphism 
\[
H^Q_2(\mathcal{P}; \mathbb{Z}[\mathcal{P}]) \longrightarrow \hatBC
\]
along with the work of Dupont and Zickert \cite{dupont-zickert} 
by using the map $\varphi: H^Q_2(\mathcal{P}; \mathbb{Z}[\mathcal{P}]) \to H^{\Delta}_3(\mathcal{P})$.

\subsection{}
\label{subsec:C_3_to_hatBC}
Let $C_{n}( \Cpm )$ be the free abelian group generated by $(n+1)$-tuples of elements of $\Cpm$.
Define a map $\partial : C_{n}(\Cpm) \rightarrow C_{n-1}(\Cpm)$ by
\[
 \partial (v_{0}, v_{1}, \cdots, v_{n}) = \sum_{i = 0}^{n} (-1)^{i} (v_{0}, v_{1}, \cdots, \widehat{v_{i}}, \cdots, v_{n}).
\]
Then $\partial$ satisfies $\partial \circ \partial = 0$.
The group $G$ acts on $C_{n}(\Cpm)$ from the left by $g (v_{0}, v_{1}, \cdots, v_{n}) = (g v_{0}, g v_{1}, \cdots, g v_{n})$.
Therefore, $C_{n}(\Cpm)$ is a left $\mathbb{Z}[G]$-module.
Let
\[
 C_{n}(\Cpm)_{G} = \mathbb{Z} \otimes_{\mathbb{Z}[G]} C_{n}(\Cpm).
\]
We denote the $n$-th homology group of $C_{\ast}(\Cpm)_{G}$ by $H_{n}(\Cpm)$.
Since $\mathcal{P}$ is identified with $\Cpm$, $C^{\Delta}_n(\mathcal{P})$ is isomorphic to $C_n(\Cpm)$.
Moreover, since $G_{\mathcal{P}}$ acts on $\mathcal{P}$ by (\ref{eq:binary_operation_on_C^2}) and the map $G_{\mathcal{P}} \to \PSLC$ is surjective, 
we can identify $C^{\Delta}_n(\mathcal{P})_{G_{\mathcal{P}}}$ with $C_n(\Cpm)_G$.
Therefore $H^{\Delta}_n(\mathcal{P})$ is isomorphic to $H_n(\Cpm)$.

Let $h : \Cpm \rightarrow \mathbb{C}P^{1}$ be the natural map defined by
\[
 \begin{pmatrix} \alpha \\ \beta \end{pmatrix} \mapsto \frac{\alpha}{\beta}.
\]
Let $ C^{h \neq}_n(\Cpm)$ be a subcomplex of $C_n( \Cpm )$ generated by $(v_0,\dots,v_n)$ satisfying $h(v_i) \neq h(v_j)$ for $i \neq j$.
We let 
\[
C_{n}^{h \neq}(\Cpm)_{G} = \mathbb{Z} \otimes_{\mathbb{Z}[G]} C_{n}^{h \neq}(\Cpm)
\]
and $H_{n}^{h \neq}(\Cpm)$ be the $n$-th homology group of $C_{\ast}^{h \neq}(\Cpm)_{G}$.

For each $v_{i} = \begin{pmatrix} \alpha_{i} \\ \beta_{i} \end{pmatrix}, v_{j} = \begin{pmatrix} \alpha_{j} \\ \beta_{j} \end{pmatrix} \in \Cpm$ 
satisfying $h(v_{i}) \neq h(v_{j})$, 
\[
\det(v_{i}, v_{j}) = \det \begin{pmatrix} \alpha_{i} & \alpha_{j} \\ \beta_{i} & \beta_{j} \end{pmatrix}
\]
is non-zero and well-defined up to sign.
We fix a sign of $\det(v_i,v_j)$ once and for all, for example, to satisfy $0 \leq \arg(\det(v_{i}, v_{j})) < \pi$. 
Then $\Log (\det (v_i,v_j))$ is well-defined and satisfies
\begin{equation}
\label{eq:invariance_under_PSLC}
\Log (\det(gv_i,gv_j)) = \Log (\det(v_i,v_j))
\end{equation}
for any $g \in \PSLC$.
Let $(v_0,v_1,v_2,v_3)$ be a generator of $C_3^{h \neq}( (\mathbb{C}^2 \setminus \{0\})/\pm )$.
Since we have
\[
\begin{split}
[h(v_0):h(v_1):h(v_2):h(v_3)] 
& = \frac{\alpha_0/\beta_0-\alpha_3/\beta_3}{\alpha_1/\beta_1-\alpha_3/\beta_3} 
  \frac{\alpha_1/\beta_1-\alpha_2/\beta_2}{\alpha_0/\beta_0-\alpha_2/\beta_2} \\
& =\pm \frac{\det({v_0},{v_3})\det({v_1},{v_2})}{\det({v_1},{v_3})\det({v_0},{v_2})}, 
\end{split}
\]
therefore 
\begin{equation}
\label{eq:log_of_crossratio_1}
\begin{split}
\Log \det({v_0},{v_3}) + \Log \det({v_1},{v_2}) 
- \Log \det({v_1},{v_3})  - \Log \det({v_0},{v_2}) \\
= \Log( [h(v_0):h(v_1):h(v_2):h(v_3)]) + p \pi i 
\end{split}
\end{equation}
for some integer $p$. Similarly, we have 
\begin{equation}
\label{eq:log_of_crossratio_2}
\begin{split}
\Log \det(v_0 , v_2) + \Log \det(v_1, v_3) 
- \Log \det(v_0,v_1)  - \Log \det(v_2, v_3) \\
= \Log([h(v_1):h(v_2):h(v_0):h(v_3)]) + q \pi i 
\end{split}
\end{equation}
for some $q \in \mathbb{Z}$.
Define
\begin{equation}
\label{C^2-to-extended}
\begin{split}
w_0 &= \Log \det(v_0 , v_3 ) + \Log \det(v_1 , v_2 )
- \Log \det(v_0 , v_2 ) - \Log \det(v_1 , v_3 ), \\
w_1 &= \Log \det(v_0 , v_2 ) + \Log \det(v_1 , v_3 ) 
- \Log \det(v_0 , v_1 ) - \Log \det(v_2 , v_3 ), \\
w_2 &= \Log \det(v_0 , v_1 ) + \Log \det(v_2 , v_3 )
- \Log \det(v_0 , v_3 ) - \Log \det(v_1 , v_2 ). 
\end{split}
\end{equation}
By equations (\ref{eq:log_of_crossratio_1}) and (\ref{eq:log_of_crossratio_2}), 
this gives a combinatorial flattening of the ideal tetrahedron with cross ratio $[h(v_0):h(v_1):h(v_2):h(v_3)]$.
Define a map $\widehat{\sigma}: C^{h \neq}_3(\Cpm) \to \hatPC$ by $\widehat{\sigma}(v_0,v_1,v_2,v_3) = (w_0,w_1,w_2)$. 
This induces a map $C^{h \neq}_3(\Cpm)_G \to \hatPC$ by (\ref{eq:invariance_under_PSLC}). 
Since $\widehat{\sigma} (\partial (v_0, \dots, v_4))$ satisfies the lifted five term relation,
we obtain a map $\widehat{\sigma}: H^{h \neq}_{3}( \Cpm ) \to \hatPC$.
Moreover, the image of this map is in $\hatBC$. 
In fact, the map $\mu: C_*^{h \neq}(\mathbb{C}^2 \backslash \{ 0 \} / \pm)_G \to \mathbb{C} \wedge_{\mathbb{Z}} \mathbb{C}$
defined by
\[
\begin{split}
(v_0, & v_1,v_2) \mapsto  \Log\det({v_0},{v_1}) \wedge \Log\det({v_0},{v_2})  \\
& - \Log\det({v_0},{v_1}) \wedge \Log\det({v_1},{v_2})  
 + \Log\det({v_0},{v_2}) \wedge \Log\det({v_1},{v_2})
\end{split}
\]
satisfies the following commutative diagram:
\[
\xymatrix{
 C_{3}^{h \neq}(\Cpm)_{G} \ar[r]^{\hskip 2.8em \widehat{\sigma}} \ar[d]^{\partial} & \hatPC \ar[d]^{\nu} \\
 C_{2}^{h \neq}(\Cpm)_{G} \ar[r]^{\hskip 2.8em \mu} & \mathbb{C} \wedge_{\mathbb{Z}} \mathbb{C}
}
\]
Therefore, we obtain a map 
\begin{equation}
\widehat{\sigma}:H^{h \neq}_{3}( \Cpm ) \longrightarrow \hatBC.
\end{equation}
Using the \emph{cycle relation} discussed in \cite{neumann04}, we can show the following.
\begin{lemma}
\label{lem:independence_of_the_sings}
$\widehat{\sigma} : H_{3}^{h \neq}((\mathbb{C}^{2} \setminus \{ 0 \}) / \pm) \rightarrow \hatBC$ does not depend on the choice of the sign of $\det(v_{i}, v_{j})$.
\end{lemma}
\begin{proof}
Fix $v, w \in (\mathbb{C}^{2} \setminus \{ 0 \}) / \pm$ satisfying $h(v) \neq h(w)$.
Let ${\det}^{\prime}$ be the determinant function $(\Cpm)^2 \to \mathbb{C}$ with another choice of sign defined by
\[
 {\det}^{\prime}(v_{i}, v_{j}) = 
 \begin{cases}
  - \det(v_{i}, v_{j}) & \mathrm{if} \ \{ v_{i}, v_{j} \} = \{ v, w \}, \\
  \det(v_{i}, v_{j}) & \mathrm{otherwise},
 \end{cases}
\]
for each $v_{i}, v_{j} \in \mathbb{C}^{2} \setminus \{ 0 \} / \pm$.
Suppose ${\widehat{\sigma}}^{\prime} : H_{3}^{h \neq}((\mathbb{C}^{2} \setminus \{ 0 \}) / \pm) \rightarrow \hatBC$ is the map associated with ${\det}^{\prime}$.
It is sufficient to show that $\widehat{\sigma}$ and ${\widehat{\sigma}}^{\prime}$ are the same.

Suppose $C$ is a cycle of $C_{3}^{h \neq}((\mathbb{C}^{2} \setminus \{ 0 \}) / \pm)$.
Let $\Delta_{1}, \cdots, \Delta_{n}$ be the simplices of $C$ one of whose edges is the edge $[v w]$ spanned by $v$ and $w$.
By the assumption that $C$ is a cycle, replacing the indices of $\Delta_{i}$ if necessary, we may assume that $\Delta_{i}$ is glued to $\Delta_{i+1}$ along one of the two faces of $\Delta_{i}$ incident to $[v w]$ for each index $i$ modulo $n$.
We label the two edges other than $[v w]$ of the common face of $\Delta_{i}$ and $\Delta_{i+1}$ as $T_{i}$ and $B_{i}$ 
(for ``top'' and ``bottom'') in such a way that $T_i$ (resp. $B_i$) incident to $v$ (resp. $w$).

By definition, the difference between $\Log \hskip 0.1em {\det}^{\prime}(v, w)$ and $\Log \det(v, w)$ is $\pm \pi i$.
Therefore, by (\ref{C^2-to-extended}), 
\begin{itemize}
 \item the log-parameters of ${\widehat{\sigma}}^{\prime}(\Delta_{i})$ at $B_{i}$ and its opposite edge in $\Delta_i$  
   are $\mp \pi i$ of those of $\widehat{\sigma}(\Delta_{i})$, 
 \item the log-parameters of ${\widehat{\sigma}}^{\prime}(\Delta_{i})$ at $T_{i}$ and its opposite edge in $\Delta_i$  
   is $\pm \pi i$ of those of $\widehat{\sigma}(\Delta_{i})$, 
 \item the log-parameter of ${\widehat{\sigma}}^{\prime}(\Delta_{i})$ at the edge $[v w]$ is equal to that of $\widehat{\sigma}(\Delta_{i})$.
\end{itemize}
Further, the sum of the log-parameters around $[v w]$ of $\widehat{\sigma}(\Delta_{i})$ is zero.
In this case, according to Lemma 6.1 of \cite{neumann04}, ${\widehat{\sigma}}^{\prime}(C)$ is equal to ${\widehat{\sigma}}^{\prime}(C)$ in $\hatP(\mathbb{C})$.
\end{proof}

\subsection{}
In Section \ref{sec_main}, we have constructed a map $\varphi : H^{Q}_{2}(\mathcal{P}; \mathbb{Z}[\mathcal{P}]) \rightarrow H^{\Delta}_{3}(\mathcal{P})$.
Thus, if we have a map from $H^{\Delta}_{3}(\mathcal{P})$, which is isomorphic to $H_3( \Cpm )$, to $\hatBC$, we obtain an element of the extended Bloch group $\hatBC$ 
associated with a shadow coloring.
We do not directly construct a map from $H_3( \Cpm )$ to $\hatBC$, instead we prove Proposition \ref{prop:lift}.
Before stating it, we need a simple observation.

\begin{lemma}
The homomorphism $H^{h \neq }_{n}( \Cpm ) \to H_{n}(\Cpm)$ induced from the inclusion is injective for $n \geq 1$.
\end{lemma}
\begin{proof}
For any $x_0, x_1 \in \Cpm$, we fix $y(x_0,x_1) \in \Cpm$ so that $h(y(x_0,x_1))$ is different from $h(x_0)$ and $h(x_1)$.
Since the map $h$ is $G$-equivariant, we can choose $y(x_0,x_1)$ to satisfy $y(g x_0, g x_1) = g y(x_0, x_1)$ for $g \in G$.
For any triple $x_0, x_1, x_2 \in \Cpm$, we fix $y(x_0,x_1,x_2) \in \Cpm$ so that $h(y(x_0,x_1,x_2))$ is different from $h(x_i)$ and $h(y(x_i,x_j))$ for all $i,j$.
We can also assume that $y$ is $G$-equivariant.
Inductively, fix $y(x_0, \cdots, x_n) \in \Cpm$ in a $G$-equivariant way so that $h(y(x_0, \cdots, x_n))$ is different from $h(y(x_{i_0}, \cdots, x_{i_k}))$ 
for any subset $\{x_{i_0}, \dots, x_{i_k}\}$ of $\{x_{0},\dots x_{n}\}$.
We can assume that $y(x_{\sigma(0)}, \cdots, x_{\sigma(n)}) = y(x_0, \cdots, x_n)$ for any permutation $\sigma$ of $\{0, \dots n\}$.

By abuse of notation, we denote $y(x_0, \cdots, x_n)$ by $x_{01 \cdots n}$ for an $n$-tuple of points $x_0, \cdots, x_n$.
We define a barycentric subdivision map $b_n: C_n(\Cpm) \to C^{h \neq}_n(\Cpm)$ by
\[
b_n( (x_0, \cdots, x_n) ) =\sum_{\sigma \in \mathfrak{S}_{n+1}} \mathrm{sgn}(\sigma) (x_{\sigma(0)}, x_{\sigma(0) \sigma(1)}, \cdots, x_{\sigma(0) \cdots \sigma(n)}).
\]
This is a chain map of $\mathbb{Z}[G]$-module chain complexes. We remark that $C^{h \neq}_0(\Cpm) = C_0(\Cpm)$ and $b_0$ is the identity map.
Let $K$ be the kernel of the augmentation map $C^{h \neq}_0(\Cpm) \to \mathbb{Z}$. 
Since the augmented chain complex of $(C^{h \neq}_n(\Cpm), \partial_n)$ is acyclic, $K \cong \mathrm{Im}( \partial_1 ) \cong \mathrm{Coker}(\partial_2)$.
Since $C^{h \neq}_n(\Cpm)$ is a free $\PSLC$-module for $n \geq 1$, 
\[
\cdots \to C^{h \neq}_2(\Cpm) \to C^{h \neq}_1(\Cpm) \to K \to 0
\]
is a free resolution of $K$.
The composition of the inclusion $C^{h \neq}_n(\Cpm) \to C_n(\Cpm)$ and the barycentric subdivision map is the identity at $K$,
the induced map $H^{h \neq }_{n}( \Cpm ) \to H_{n}(\Cpm) \to H^{h \neq }_{n}( \Cpm )$ is an isomorphism for $n \geq 1$.
Thus $H^{h \neq }_{n}( \Cpm ) \to H_{n}(\Cpm)$ is an injection for $n \geq 1$.
\end{proof}

Thus we can ask whether or not the image of $\varphi : H^{Q}_{2}(\mathcal{P}; \mathbb{Z}[\mathcal{P}]) \rightarrow H^{\Delta}_{3}(\mathcal{P})$ is in $H^{h \neq}_3( \Cpm )$.
\begin{proposition}
\label{prop:lift}
Let $C$ be a cycle of $C^{Q}_{2}(\mathcal{P}; \mathbb{Z}[\mathcal{P}])$.
Then the image $[\varphi(C)]$ is in $H^{h \neq}_3 (\Cpm)$.
As a result we obtain a homomorphism 
\[
\psi : H^{Q}_{2}(\mathcal{P}; \mathbb{Z}[\mathcal{P}]) \longrightarrow H_{3}^{h \neq}(\Cpm).
\]
\end{proposition}

\begin{proof}
Let $C$ be a cycle of $C^{Q}_{2}(\mathcal{P}; \mathbb{Z}[\mathcal{P}])$.
Then there is a shadow coloring $\mathcal{S}$ of a link diagram $D$ on an orientable surface satisfying $C(\mathcal{S}) = C$ (Remark \ref{rem:general_cycles}).
We assume that $C(S) = \sum_{c} \varepsilon_{c} r_{c} \otimes (x_{c}, y_{c})$.
Here, $c$ represents a crossing of $D$, $\varepsilon_{c}$ is $1$ or $-1$ depending on whether $c$ is positive or negative respectively, 
and $x_{c}, y_{c}$ and $r_{c}$ denote colors around $c$ with respect to $\mathcal{S}$ as illustrated in Figure \ref{fig:cycle_from_shadow_coloring}.
By (\ref{eq:does_not_degenerate}), we have
\begin{eqnarray}
 \varphi \left( \left[ \sum_{c} \varepsilon_{c} r_{c} \otimes (x_{c}, y_{c}) \right] \right) & & \nonumber \\
 &\hskip -15em = & \hskip -7em \Bigg[ \sum_{c} \varepsilon_{c}( (p, r_{c}, x_{c}, y_{c}) - (p, r_{c} \ast x_{c}, x_{c}, y_{c}) \nonumber \\
\label{eq:retriangulation}
 &\hskip -15em  & \hskip -7em \hskip 4em - \> (p, r_{c} \ast y_{c}, x_{c} \ast y_{c}, y_{c}) + (p, r_{c} \ast (x_{c} y_{c}), x_{c} \ast y_{c}, y_{c})) \Bigg] \\
 &\hskip -15em  = & \hskip -7em \Bigg[ \sum_{c} \varepsilon_{c}( (p, r_{c} \ast x_{c}, r_{c}, x_{c}) - (p, r_{c} \ast x_{c}, r_{c}, y_{c}) \nonumber \\
 &\hskip -15em  & \hskip -7em \hskip 4em - \> (p, r_{c} \ast (x_{c} y_{c}), r_{c} \ast y_{c}, x_{c} \ast y_{c}) 
 + (p, r_{c} \ast (x_{c} y_{c}), r_{c} \ast y_{c}, y_{c})) \Bigg]. \nonumber
\end{eqnarray}
Recall that the homology class $[C(\mathcal{S})] \in H^{Q}_{2}(\mathcal{P}; \mathbb{Z}[\mathcal{P}])$ does not depend on the choice of a region coloring 
(Theorem \ref{thm:region_and_conjugacy_invariance} and Remark \ref{rem:general_cycles}), 
and the homomorphism $\varphi$ does not depend on the choice of $p \in \mathcal{P}$ (Proposition \ref{prop:does_not_depend}).
We can change the region coloring $\{r_c\}_{c} \subset \mathcal{P}$ and the point $p \in \mathcal{P}$ 
so that the simplices appeared in the fourth and fifth lines of (\ref{eq:retriangulation}) do not degenerate as follows.
Recall that a region coloring is uniquely determined by the color, say $r$, of one region.
If we replace $r$ so that $h(r)$ is away from the fixed points corresponding to the meridians, 
then $\{h(r_c)\}_{c}$ does not intersect the fixed points of the arc colors since each $h(r_c)$ lies in the orbit of $h(r)$ under the action of the arc colors.
Now we have $h(r_c) \neq  h(r_c * x_c)$.
In fact, if $h(r_c) = h(r_c * x_c)$, then $h(r_c) = x_c h(r_c)$ by (\ref{eq:binary_operation_on_C^2}).
Thus $h(r_c)$ must be the fixed point of $x_c$, which implies $h(r_c) = h(x_c)$.
This contradicts the choice of $r$. 
Hence $(h(r_c*x_c), h(r_c), h(x_c))$ and $(h(r_{c} \ast x_{c}), h(r_{c}), h(y_{c}))$ are distinct triples of points in $\mathbb{C}P^1$.
We can also show that 
$(h(r_{c} \ast (x_{c} y_{c}) ), h(r_{c} \ast y_{c}), h(x_{c} \ast y_{c}))$ and  $(h(r_{c} \ast (x_{c} y_{c}), h(r_{c} \ast y_{c}), h(y_{c}))$ are
distinct triples.
If we further choose $p$ so that $h(p)$ is away from these points, then 
\[
\begin{split}
 &(h(p), h(r_{c} \ast x_{c}), h(r_{c}), h(x_{c})), \enskip (h(p), h(r_{c} \ast x_{c}), h(r_{c}), h(y_{c})), \\
 &(h(p), h(r_{c} \ast (x_{c} y_{c}) ), h(r_{c} \ast y_{c}), h(x_{c} \ast y_{c})), \enskip (h(p), h(r_{c} \ast (x_{c} y_{c})), h(r_{c} \ast y_{c}), h(y_{c})) \\
\end{split}
\]
are non-degenerate tetrahedra for all $c$.
\end{proof}

\begin{remark}
\label{rem:on_the_tirangulation}
In the proof of Proposition \ref{prop:lift}, $\varphi(C(\mathcal{S}))$ may contain degenerate ideal tetrahedra in general, 
so we need to change $\varphi(C(\mathcal{S}))$ by (\ref{eq:retriangulation}).
(In this paper, we mean that an ideal tetrahedron is \emph{degenerate} if two of their ideal vertices coincide.
In this sense, a flat tetrahedron is non-degenerate unless the associated cross ratio is equal to $0$ or $1$.)
For example, if we give an arc coloring $\mathcal{A}$ by one element of $\mathcal{P}$ (in this case $\rho_{\mathcal{A}}$ is a reducible representation), 
then we can not separate the points $h(x_c)$ and $h(y_c)$ whatever the region coloring is.
But we can choose an appropriate region coloring to separate $h(r_c)$, $h(r_c * x_c)$ and $h(x_c)$ as in the proof.
We will see in \S \ref{subsec:proof_of_theorem} that the change given by (\ref{eq:retriangulation}) corresponds to 
a retriangulation of the link complement (see Figure \ref{fig:retriangulation}).
The point is that we can change the colors corresponding to the bottom vertex in Figure \ref{fig:retriangulation} not freely but in a well-controlled way. 
We also remark that the modified triangulation is still not an ideal triangulation in the usual sense.
It has two non-ideal points, which are called the north and south poles in \cite{weeks}. 
Thus Proposition \ref{prop:lift} does not give a non-degenerate ideal triangulation in the usual sense.
\end{remark}

By Proposition \ref{prop:lift}, we obtain a cohomology class
\[
 [\mathrm{cvol}] \in H_{Q}^{2}(\mathcal{P}; \mathrm{Hom}(\mathbb{Z}[\mathcal{P}], \mathbb{C} / \pi^{2} \mathbb{Z}))
\]
defined by $\langle \mathrm{cvol}, C \rangle = R(\widehat{\sigma}(\psi([C])))$ for $C \in C^{Q}_{2}(\mathcal{P}; \mathbb{Z}[\mathcal{P}])$.
Here $R : \hatBC \rightarrow \mathbb{C} / \pi^{2} \mathbb{Z}$ is the map defined in Section \ref{sec_extended_bloch}.
This cohomology class $[\mathrm{cvol}]$ gives the complex volume of a hyperbolic link as follows.

Let $L$ be an oriented link in $S^3$. 
Let $D$ be a diagram of $L$ and $\mathcal{S} = (\mathcal{A}, \mathcal{R})$ a shadow coloring of $D$ with respect to $\mathcal{P}$.
Recall that, associated with a shadow coloring $\mathcal{S}$, 
we have a parabolic representation $\xi \circ \rho_{\mathcal{A}} : \pi_{1}(S^{3} \setminus L) \rightarrow \PSLC$ where $\xi : G_{\mathcal{P}} \rightarrow \PSLC$
be the surjective homomorphism induced by the natural inclusion $\mathcal{P} \rightarrow \PSLC$.
On the other hand, the shadow coloring $\mathcal{S}$ determines a cycle $[C(\mathcal{S})]$ in $H^Q_2(\mathcal{P}, \mathbb{Z}[\mathcal{P}])$.

\begin{theorem}
\label{thm:volume_and_chern_simons}
Suppose $L$ is an oriented link in $S^{3}$ and $\mathcal{S} = (\mathcal{A}, \mathcal{R})$ a shadow coloring of a diagram of $L$ with respect to $\mathcal{P}$.
Then $\widehat{\sigma}(\psi( [C(\mathcal{S})] ) ) \in \hatBC$ is equal to the invariant defined by Neumann.
If $L$ is hyperbolic and $\mathcal{A}$ corresponds to the discrete faithful representation $\xi \circ \rho_{\mathcal{A}} : \pi_{1}(S^{3} \setminus L) \rightarrow \PSLC$, then
\[
 \langle [\mathrm{cvol}], [C(\mathcal{S})] \rangle = i \hskip 0.1em (\mathrm{Vol}(S^{3} \setminus L) + i \hskip 0.1em \mathrm{CS}(S^{3} \setminus L)),
\]
where $\mathrm{CS}(S^{3} \setminus L)$ is Meyerhoff's extension of the Chern-Simons invariant to cusped hyperbolic manifolds in \cite{meyerhoff}.
\end{theorem}

Theorem \ref{thm:volume_and_chern_simons} enables us to compute the complex volume of a hyperbolic link only from a link diagram 
(see Section \ref{sec_example} for actual computations).
To prove Theorem \ref{thm:volume_and_chern_simons}, we review the work of Neumann \cite{neumann04} in the next subsection.

\begin{remark}
Recall that the fundamental class $[C(\mathcal{S})] \in H^{Q}_{2}(\mathcal{P}; \mathbb{Z}[\mathcal{P}])$ derived from a shadow coloring $\mathcal{S}$ is 
determined by the conjugacy class of $\xi \circ \rho_{\mathcal{A}}$.
Therefore, $\widehat{\sigma}(\psi( [C(\mathcal{S})] ) ) \in \hatBC$ is clearly an invariant of oriented links with parabolic representations.
If we reverse the orientation of a component of the link, 
there exists an arc coloring which induces the same representation $\xi \circ \rho_{\mathcal{A}}: \pi_1(S^3 \setminus L) \to \PSLC$ 
by assigning to each arc $\alpha_i$ of the component the inverse of $\mathcal{A}(\alpha_i)$. 
We do not know whether $[C(\mathcal{S})] \in H^{Q}_{2}(\mathcal{P}; \mathbb{Z}[\mathcal{P}])$ depends on the orientation of the link. 
But $\widehat{\sigma}(\psi( [C(\mathcal{S})] ) ) \in \hatBC$ does not depend on the orientation of the link.
\end{remark}

\subsection{}
\label{subsec:flattening}
Let $P$ be the subgroup of $G$ consisting of upper triangular matrices with $1$ on the diagonal.
Let $BG$ and $BP$ be the classifying spaces of $G$ and $P$ as discrete groups respectively. 
Neumann showed that the long exact sequence for the pair $(BG,BP)$ is simplified to the short exact sequence 
\[
 0 \longrightarrow H_{3}(BG; \mathbb{Z}) \longrightarrow H_{3}(BG, BP; \mathbb{Z}) \longrightarrow H_{2}(BP; \mathbb{Z}) \longrightarrow 0, 
\]
and there is a splitting
\[
 s : H_{3}(BG, BP; \mathbb{Z}) \longrightarrow H_{3}(BG; \mathbb{Z}).
\]

Let $M$ be a compact oriented 3-manifold and $\rho$ be a representation of $\pi_1(M)$ into $\PSLC$ 
which sends each peripheral subgroup to a parabolic subgroup, i.e. a conjugate of $P$. 
For each peripheral subgroup $H_i$, take an element $g_i$ satisfying $g_i^{-1} \rho(H_i) g_i \subset P$. 
Zickert called a choice of such elements a \emph{decoration} of $\rho$
and showed that a pair of $\rho$ and a decoration determines a fundamental class $F$ in $H_3(BG, BP;\mathbb{Z})$ \cite[Theorem 5.13]{zickert}.

On the other hand, there is an ideal triangulation of $M$ which induces the representation $\rho$ by a developing map.
Neumann defined in \cite{neumann04} (see also \cite{zickert}) an invariant $\widehat{\beta}(M) \in \hatBC$
for such an ideal triangulation and showed that it only depends on $M$ and $\rho$. 
This $\widehat{\beta}(M)$ coincides with $s(F)$ under the isomorphism $\lambda$.
(Thus $s(F) \in H_3(BG;\mathbb{Z})$ does not depend on the choice of the decoration.) 
Theorem \ref{thm:volume_and_chern_simons} means that $\widehat{\sigma}(\psi( [C(\mathcal{S})] ) )$ is equal to $\widehat{\beta}(M)$.

Let $L$ be a hyperbolic link in $S^3$ and $N(L)$ a regular neighborhood of $L$. 
For $M = S^3 \setminus N(L)$ and the discrete faithful representation $\rho : \pi_1(M) \to \PSLC$,
$R(\lambda(s(\rho_*([M, \partial M]))))$ is equal to $i(\mathrm{Vol}(S^3 \setminus L)+ i\mathrm{CS}(S^3 \setminus L))$ \cite[Corollary 14.6]{neumann04}.

We will recall the definition of $\widehat{\beta}$ in detail.
Let $K$ be a CW-complex obtained by gluing 3-simplices along their faces.
Let $K^{(i)}$ be the $i$-skeleton of $K$.
Suppose $N(K^{(0)})$ is a regular neighborhood of $K^{(0)}$.
By definition, $K \setminus N(K^{(0)})$ is a compact 3-manifold with boundary.
We denote it by $M$. 
The complex $K$ is called an ideal triangulation of $M$. 
We assume that the manifold $M$ has an orientation.
We further assume that each 3-simplex has an ordering of the vertices so that these orderings agree on common faces. 
The ordering of vertices enables us to define a combinatorial flattening for each 3-simplex of $K$.
Let $\Delta_i$ be the 3-simplices of $K$.
If the vertex ordering of $\Delta_i$ coincides with the orientation of the manifold $M$,
let $\varepsilon_i = 1$, if not, $\varepsilon_i = -1$.

Let $z_i \in \mathbb{C} \setminus \{0,1\}$ be the complex parameter of $\Delta_i$.
We assume that $z_i$ satisfies the gluing condition, i.e. the product of complex parameters around any 1-simplex of $K$ is $1$.
Then the developing map constructed from $z_i$ induces a $\PSLC$-representation of $\pi_1(M)$.
Assign a combinatorial flattening  $[z_i;,p_i,q_i]$ to each $\Delta_i$.
The log-parameter of an edge $E$ of $\Delta_i$ has a form $\Log(w) + s \pi i$ with $w \in \mathbb{C} \setminus \{0, 1\}$ and some integer $s$.
We denote $\delta_E \equiv s \mod 2$, 
in other words,
\[
 \delta_{E} \equiv 
 \begin{cases}
  p_{i} \mod 2 & \textrm{if the edge $E$ corresponds to $z_{i}$}, \\
  q_{i} \mod 2 & \textrm{if the edge $E$ corresponds to $\frac{1}{1 - z_{i}}$}, \\
  p_{i} + q_{i} + 1 \mod 2 & \textrm{if the edge $E$ corresponds to $1 - \frac{1}{z_{i}}$}.
 \end{cases}
\]
for $[z_i; p_i, q_i]$.

Let $\gamma$ be a closed path in $K$.
We call $\gamma$ \emph{normal} if $\gamma$ does not meet any 0-simplex or 1-simplex and intersects with each 2-simplex transversely.
We can deform any path to be normal.
When a normal path $\gamma$ through a $3$-simplex, entering and departing at different faces, 
there is a unique edge $E$ of the $3$-simplex between the faces.
We say that $\gamma$ \emph{passes} this edge $E$.

Let $\gamma$ be a normal path in $K$.
The \emph{parity} along $\gamma$ is the sum
\[
\sum_E \delta_E \mod 2
\]
of the parities of all the edges $E$ that $\gamma$ passes.

Let $\gamma$ be a normal path on a small neighborhood of a 0-simplex of $K$.
Let $i(E)$ be the index $i$ of the simplex to which the edge $E$ belongs.
The \emph{log-parameter} along $\gamma$ is the sum
\[
\sum_E \pm \varepsilon_{i(E)} l_E,
\]
where $E$ runs all the edges that $\gamma$ passes through and the sign $\pm$ is $+$ or $-$ according as the edge $E$ is passed in a counterclockwise 
or clockwise as viewed from the vertex.

\begin{definition}[\cite{neumann04}, \cite{zickert}]
A flattening $[z_i;p_i,q_i]$ of $K$ is called  
\begin{itemize}
\item \emph{strong flattening}, 
if the log-parameter along any normal path on a neighborhood of any 0-simplex is zero 
and the parity along any normal path in $K$ is zero; 
\item \emph{semi-strong flattening}, 
if the log-parameter along any normal path on a neighborhood of any 0-simplex is zero.
\end{itemize}
\end{definition}

For a strong flattening $[z_i;p_i,q_i]$, Neumann defined an element
\[
\hatbeta(M) = \sum_i \varepsilon_i  [z_i; p_i, q_i].
\]
\begin{proposition}[Lemma 10.1 of \cite{neumann04}]
\label{prop:neumann2}
The choice of strong flattening of $K$ does not affect the resulting element 
$\hatbeta(M) =\sum_i \varepsilon_i  [z_i; p_i, q_i]$.
\end{proposition}

We will use the following proposition later.
\begin{proposition}[Corollary 5.4 of \cite{neumann04}]
\label{prop:neumann3}
For a normal path in a neighborhood of a 0-simplex, 
if the flattening condition for log-parameters is satisfied, then so is the parity condition.
\end{proposition}

\subsection{Proof of Theorem \ref{thm:volume_and_chern_simons}}
\label{subsec:proof_of_theorem}
As discussed in \S \ref{subsec:flattening}, 
we only have to show that there exists an ideal triangulation of $S^3 \setminus L$ with a strong flattening 
which represents $\widehat{\sigma}(\psi([C(\mathcal{S})]))$.

\begin{figure}[htb]
\begin{center}
\includegraphics[scale=0.41]{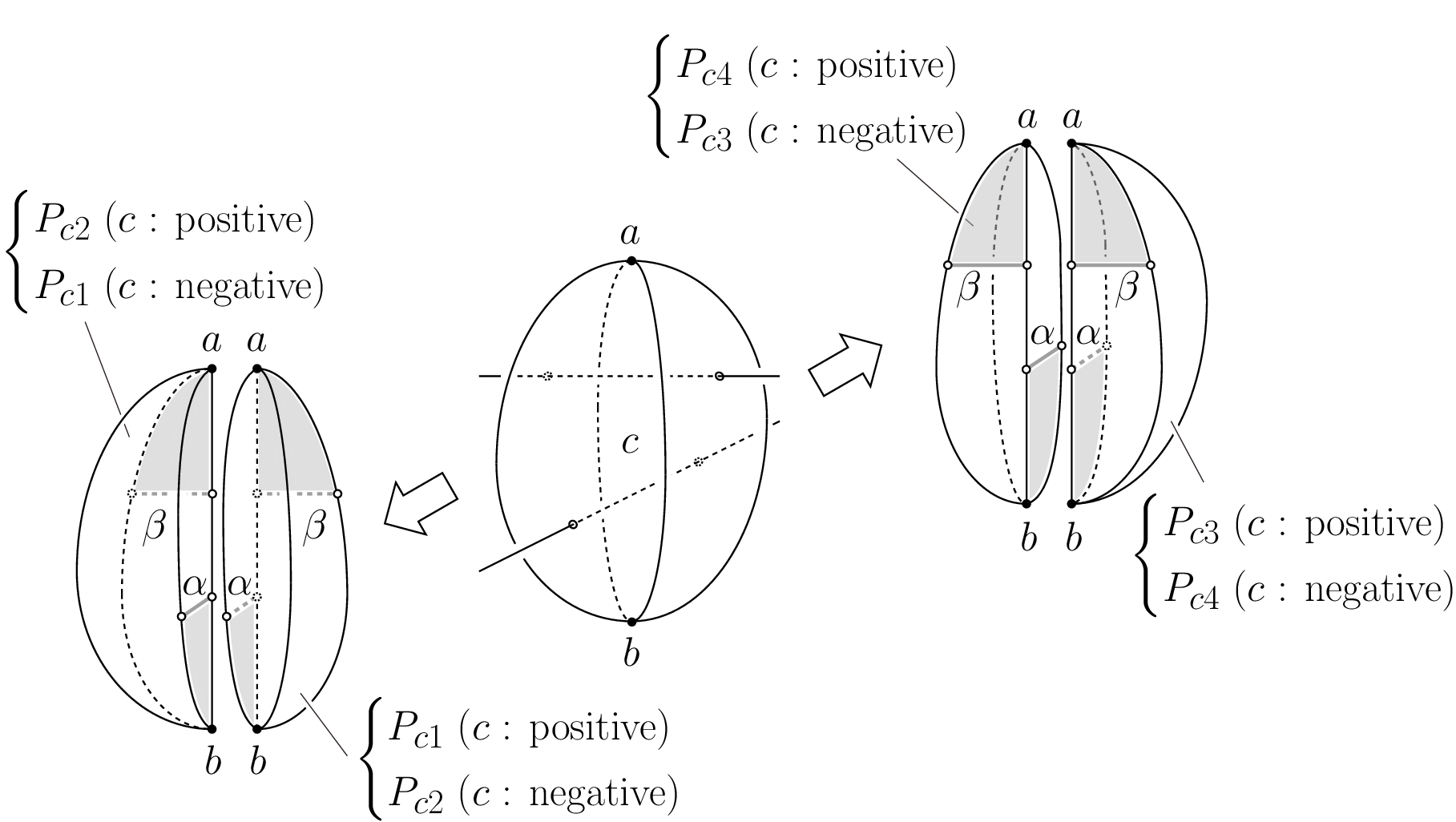}
\end{center}
\caption{Ideal tetrahedra at the crossing $c$.}
\label{fig:cut_banana}
\end{figure}

Choose and fix two points $a, b \in S^{3} \setminus L$.
Put the link diagram $D$ on a $2$-sphere which divides $S^{3}$ into two connected components containing $a$ or $b$ respectively.
Take a dual graph of $D$ on the $2$-sphere, and consider its suspension with respect to $a$ and $b$.
Then $S^{3} \setminus L$ is decomposed into thin regions,
each of which further decomposed into four pieces $P_{ci}$  as depicted in Figure \ref{fig:cut_banana}.
Compressing shaded faces into edges (then each of the edges $\alpha$ and $\beta$ degenerates into a point), 
we obtain an ideal triangulation $K$ of $S^3 \setminus L$ (see also \cite{inoue10_1}).
We remark that this decomposition of $S^{3} \setminus (L \cup \{ \mathrm{two \ solid \ balls} \})$ coincides with the one shown in \cite{weeks}.

\begin{figure}[htb]
\begin{center}
\includegraphics[scale=0.31]{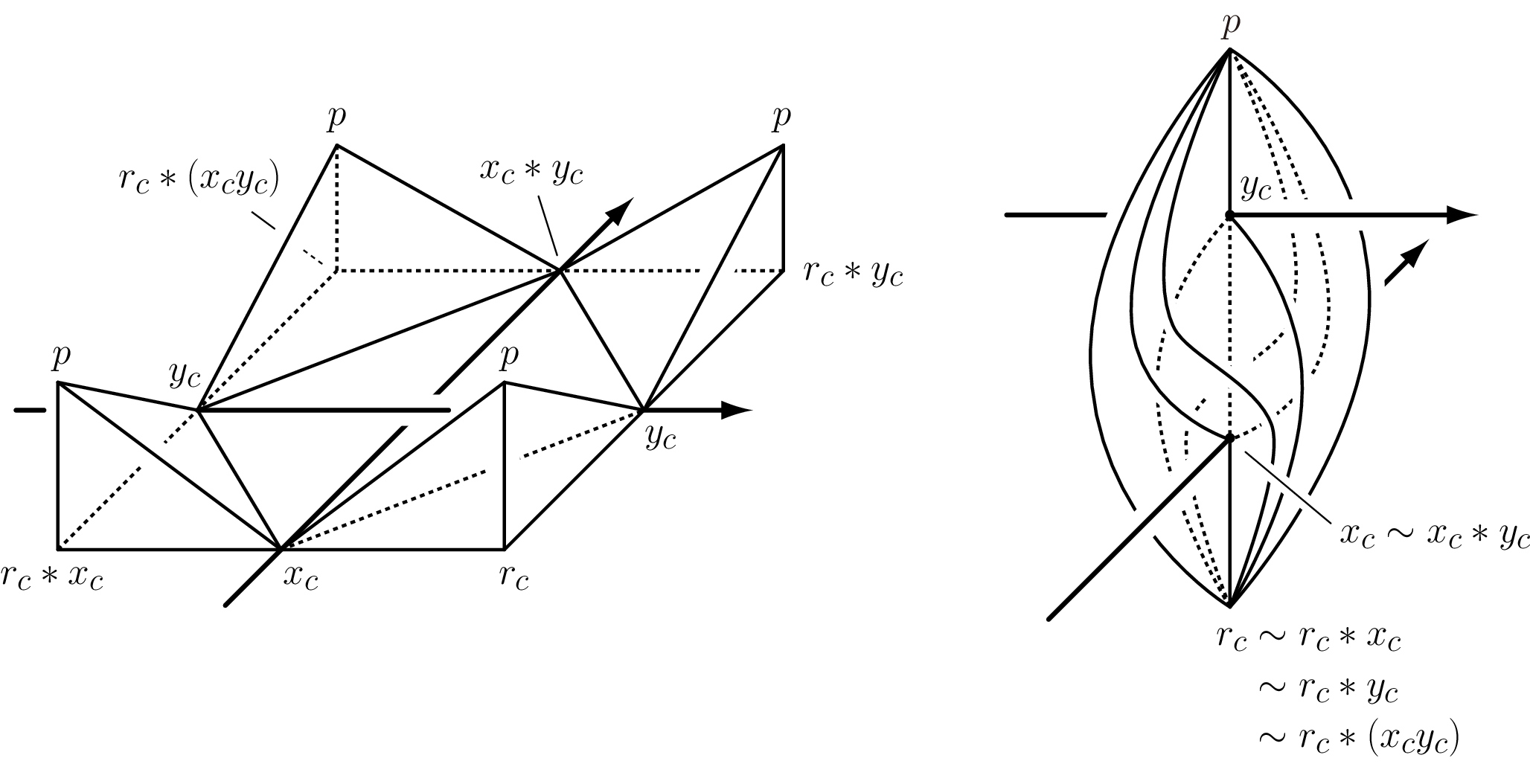}
\end{center}
\caption{The two upper faces $(p,x_c,y_c)$ (resp. $(p, x_c*y_c, y_c)$) are glued in pairs.
The lower face $(r_c,x_c,y_c)$ (resp. $(r_c*x_c,x_c,y_c)$) is glued to $(r_c*y_c, x_c*y_c, y_c)$ (resp. $(r_c*(x_cy_c), x_c*y_c,y_c)$) by the action of $y_c$.}
\label{fig:triangulation}
\end{figure}

First, we assume that $\varphi(C(\mathcal{S}))$ is in $C^{h \neq}_3(\mathbb{C}^2 \setminus \{0\}) /\pm)_G$.
Recall that $\varphi(C(\mathcal{S}))$ is given by
\begin{eqnarray*}
\begin{split}
 \sum_{c} \varepsilon_{c} & ( (p, r_{c}, x_{c}, y_{c}) - (p, r_{c} \ast x_{c}, x_{c}, y_{c}) \\
 & \enskip - (p, r_{c} \ast y_{c}, x_{c} \ast y_{c}, y_{c}) + (p, r_{c} \ast (x_{c} y_{c}), x_{c} \ast y_{c}, y_{c})).
\end{split}
\end{eqnarray*}
We let 
\[
\begin{split}
[z_{c 1}; p_{c 1}, q_{c 1}] &= \widehat{\sigma}((p, r_{c}, x_{c}, y_{c})),  \\
[z_{c 2}; p_{c 2}, q_{c 2}] &= \widehat{\sigma}((p, r_{c} \ast x_{c}, x_{c}, y_{c})), \\
[z_{c 3}; p_{c 3}, q_{c 3}] &= \widehat{\sigma}((p, r_{c} \ast y_{c}, x_{c} \ast y_{c}, y_{c})), \\
[z_{c 4}; p_{c 4}, q_{c 4}] &= \widehat{\sigma}((p, r_{c} \ast (x_{c} y_{c}), x_{c} \ast y_{c}, y_{c})). \\
\end{split}
\]
For each ideal tetrahedron corresponding to $P_{ci}$, 
assign an ordering of the vertices by $(a,b,\alpha,\beta)$ in Figure \ref{fig:cut_banana} and give the combinatorial flattening $[z_{c i}; p_{c i}, q_{c i}]$.
We can check from Figure \ref{fig:cut_banana} that the face parings between $P_{ci}$ preserve the ordering of the vertices. 
We can also check that the gluing pattern is compatible with the one given by the boundary map of $C^{h \neq}_3(\mathbb{C}^2 \setminus \{0\}) /\pm)_G$ 
(see Figure \ref{fig:triangulation}).
Thus $z_{c i}$'s satisfy the gluing condition around any 1-simplex of $K$ and gives the representation $\rho_{\mathcal{A}}$ via the developing map.
(We refer the reader to \cite{kabaya} for more details on the relation between the representation and the gluing condition.)
Therefore 
\[
 \widehat{\sigma}(\varphi([C(\mathcal{S})])) = \sum_{c} \sum_{i = 1}^{4} \varepsilon_{c i} [z_{c i}; p_{c i}, q_{c i}]
\]
and $\widehat{\sigma}(\psi([C(\mathcal{S})]))$ is a flattening of $K$.

\begin{figure}[htb]
\begin{center}
\includegraphics[scale=0.31]{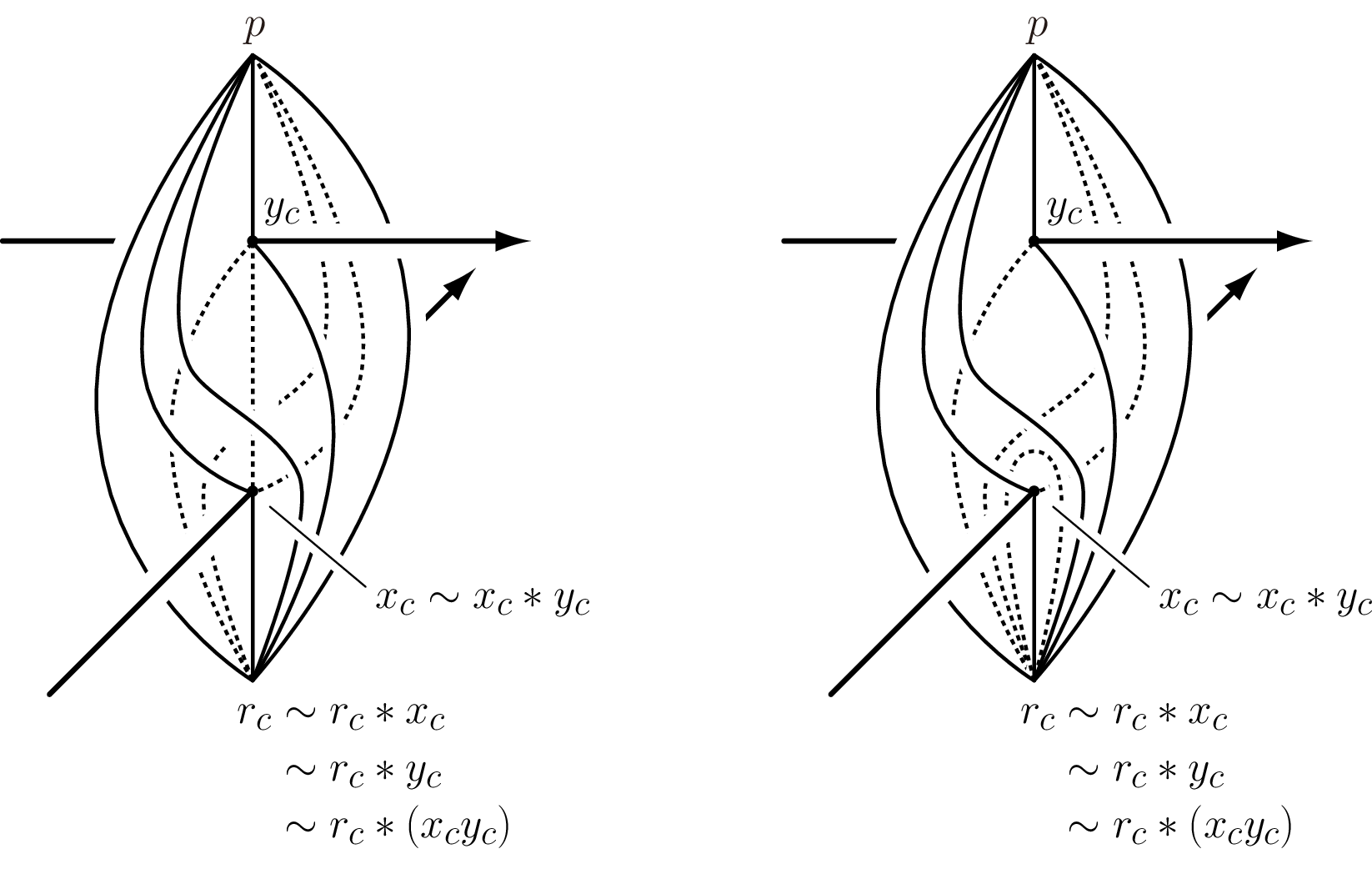}
\end{center}
\caption{The modification given by (\ref{eq:retriangulation}) corresponds to the retriangulation which exchanges the central edge diagonally.} 
\label{fig:retriangulation}
\end{figure}

If $\varphi(C(\mathcal{S}))$ is not in $C^{h \neq}_3(\mathbb{C}^2 \setminus \{0\}) /\pm)_G$, we change the triangulation as in Figure \ref{fig:retriangulation}, 
which is compatible with the operation given by (\ref{eq:retriangulation}).
Therefore we can deform $\varphi(C(\mathcal{S}))$ to be in $C^{h \neq}_3(\mathbb{C}^2 \setminus \{0\}) /\pm)_G$ 
by replacing the region coloring appropriately as discussed in the proof of Proposition \ref{prop:lift}.

We show the flattening defined above is a strong flattening.
\begin{proposition}
The flattening of $K$ is a semi-strong flattening.
\end{proposition}
\begin{proof}
We follow the proof of Theorem 6.2 of \cite{zickert}.
Let $\Delta = (v_0,v_1,v_2,v_3)$ be a simplex of $K$.
Denote $\Log(\det(v_i,v_j))$ by $c_{ij}$.
Then the log-parameter of the edge $[v_0v_1]$ is $c_{03} + c_{12} -c_{02} -c_{13}$.
In other words, the log-parameter of the edge $[v_0v_1]$ is a signed sum of $c_{ij}$ 
that do not correspond to $[v_0v_1]$ nor the opposite edge $[v_2v_3]$. 
Let $\gamma$ be a normal path on a neighborhood of a 0-simplex of $K$.
As depicted in Figure \ref{fig:flattening}, the sum of log-parameters along $\gamma$ cancels out.
\begin{figure}[htb]
\begin{center}
\includegraphics[scale=0.32]{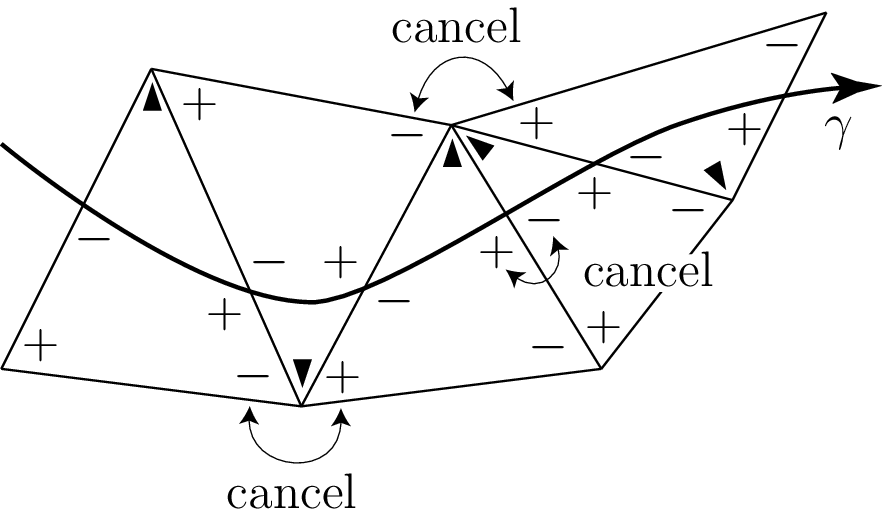}
\end{center}
\caption{Cancellation of the log-parameter along $\gamma$.}
\label{fig:flattening}
\end{figure}
\end{proof}

\begin{proposition}
\label{prop:parity}
For any normal path $\gamma$ in $K$, the parity along $\gamma$ is zero.
\end{proposition}
\begin{proof}
Let $M = K \setminus N(K^{(0)})$.
Since the sum of log-parameters around any 1-simplex of $K$ is zero, 
the parity gives an element of $H^1(M;\mathbb{Z}/2) = \mathrm{Hom}(H_1(M), \mathbb{Z}/2)$.
Now $M$ is obtained by a link complement by removing two solid balls, $H_1(M;\mathbb{Z}/2)$ is generated by normal paths on $\partial M$. 
So we only have to check that the parity along any normal path $\gamma$ on $\partial M$ is zero.
Since $K$ satisfies semi-strong flattening condition, the log-parameter along any normal path $\gamma$ on $\partial M $ is zero.
By Proposition \ref{prop:neumann3}, the parity along $\gamma$ is also zero.
\end{proof}

\subsection{$n$-dimensional  hyperbolic volume}
We end this section by discussing the $n$-dimensional hyperbolic volume.
Let $\mathbb{H}^n$ be the $n$-dimensional hyperbolic space and $\mathrm{Isom}^+(\mathbb{H}^n)$ the group of orientation preserving isometries.
Let $\mathcal{P}_n$ be the set of all parabolic elements of $\mathrm{Isom}^+(\mathbb{H}^n)$, then $\mathcal{P}_n$ is a conjugation quandle.
Since any parabolic transformation has a unique fixed point on the ideal boundary $\partial \overline{\mathbb{H}^n}$, 
we can define a map $\mathcal{P}_n \to \partial \overline{\mathbb{H}^3}$ by sending a parabolic element to its fixed point.
This induces a map $C^{\Delta}_n(\mathcal{P}_n) \to C_n(\partial \overline{\mathbb{H}^3})$, 
where $C_n(\partial \overline{\mathbb{H}^3})$ is the free abelian group generated by $(n+1)$-tuples of elements of $\partial \overline{\mathbb{H}^n}$. 
Let $\mathrm{vol}_n: C_n(\partial \overline{\mathbb{H}^3}) \to \mathbb{R}$ be the signed volume of the convex hull of the $n+1$ points.
Composting these maps, we obtain a map $\mathrm{vol}_n: C^{\Delta}_n(\mathcal{P}_n)_{G_{\mathcal{P}_n}} \to \mathbb{R}$ since the volume is invariant under isometries.
By Stokes' theorem, $\mathrm{vol}_n$ vanishes on the boundaries of $C^{\Delta}_n(\mathcal{P})$.
Therefore $n$-dimensional hyperbolic volume gives rise to an $n$-cocycle in the simplicial quandle cohomology 
and also an $(n-1)$-cocycle in the rack cohomology  $H_{R}^{n-1}(\mathcal{P}_n; \mathrm{Hom}(\mathbb{Z}[\mathcal{P}_n], \mathbb{R}))$ by Theorem \ref{thm:homomorphisms}.
Moreover, $\mathrm{vol}_n$ is a quandle cocycle, since the volume is zero if two ideal vertices coincide. 
In the case $n=3$, this is the quandle cocycle obtained in \cite{inoue10_1}.

\section{Example}
\label{sec_example}
\begin{figure}[htb]
\begin{center}
\includegraphics[scale=0.32]{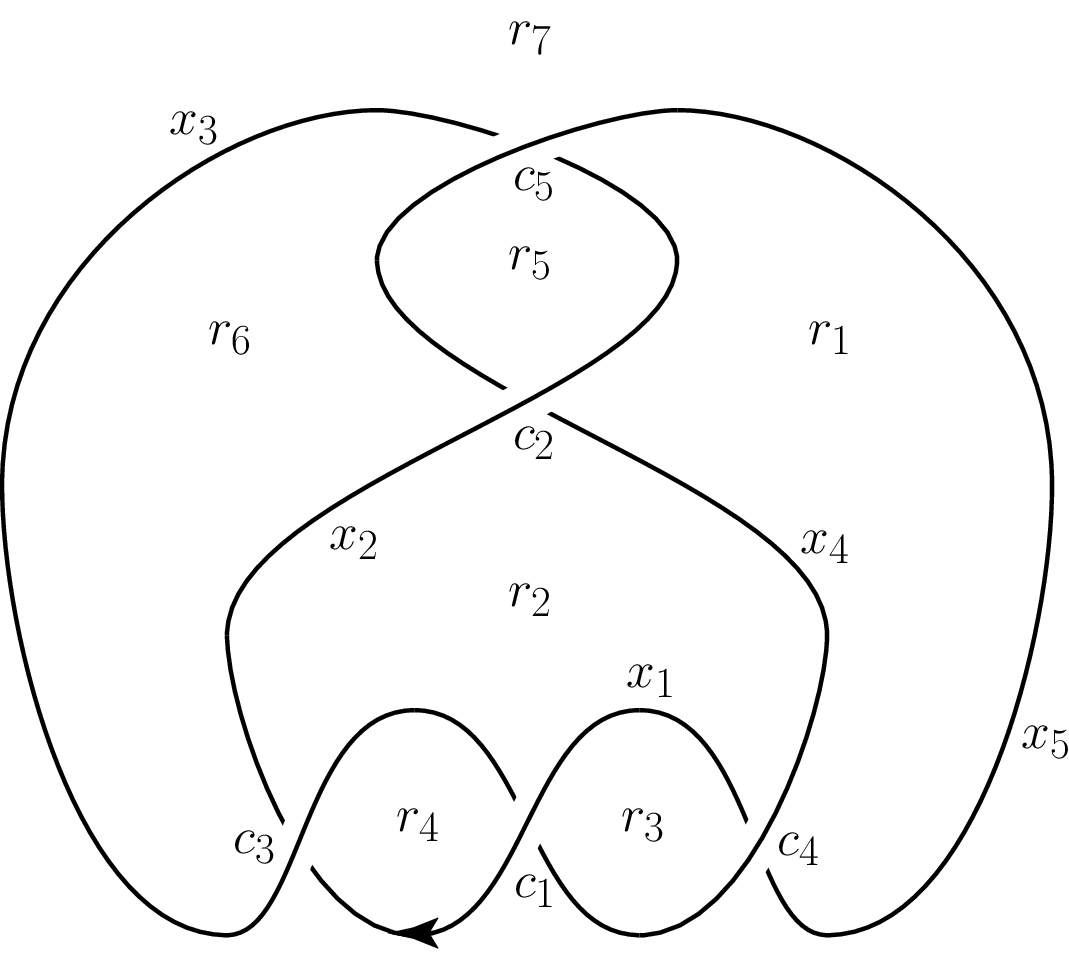}
\end{center}
\caption{A diagram $D$ of $5_2$ knot.}
\label{fig:5_2}
\end{figure}

In this section, we compute the complex volume of $5_{2}$ knot, which is hyperbolic, from a concrete shadow coloring using Theorem \ref{thm:volume_and_chern_simons}.
Define a diagram $D$ of $5_{2}$ as depicted in Figure \ref{fig:5_2}.

We first construct a shadow coloring $\mathcal{S} = (\mathcal{A}, \mathcal{R})$ of $D$ with respect to $\Cpm$ such that the parabolic representation 
$\xi \circ \rho_{\mathcal{A}} : \pi_{1}(S^{3} \setminus 5_{2}) \rightarrow \PSLC$ derived from $\mathcal{A}$ is discrete and faithful.
Suppose $x_{1}, x_{2}, \cdots, x_{5} \in \Cpm$ are colors of arcs with respect to $\mathcal{A}$ (see Figure \ref{fig:5_2}).
We let
\[
 x_{1} = \begin{pmatrix} 1 \\ 0 \end{pmatrix}, \enskip x_{4} = \begin{pmatrix} 0 \\ t \end{pmatrix}
\]
with some $t \in \mathbb{C} \setminus \{ 0 \}$, where $t$ is well-defined up to sign.
From the relations $x_{3} = x_{4} \ast^{-1} x_{1}$, $x_{2} = x_{1} \ast x_{3}$, and $x_{5} = x_{1} \ast^{-1} x_{4}$ 
at crossings $c_{1}$, $c_{3}$, and $c_{4}$ respectively (see the left-hand side of Figure \ref{fig:coloring}), we have
\[
 x_{3} = \begin{pmatrix} t \\ t \end{pmatrix}, \enskip x_{2} = \begin{pmatrix} t^{2} + 1 \\ t^{2} \end{pmatrix}, \enskip x_{5} = \begin{pmatrix} 1 \\ -t^{2} \end{pmatrix} 
\]
by (\ref{eq:binary_operation_on_C^2}) and (\ref{eq:inverse_binary_operation_on_C^2}).
Further, from  the relations $x_{5} = x_{4} \ast x_{2}$ and $x_{3} = x_{2} \ast x_{5}$ at the crossings $c_{2}$ and $c_{5}$ respectively, 
we have the relations
\[
 \begin{pmatrix} 1 \\ -t^{2} \end{pmatrix} = \begin{pmatrix} -t (1 + t^{2})^{2} \\ t (-t^{4} - t^{2} + 1) \end{pmatrix}, \enskip 
\begin{pmatrix} t \\ t \end{pmatrix} = \begin{pmatrix} -t^{4} - t^{2} + 1 \\ t^{2} (t^{2} + 1)^{2} \end{pmatrix}.
\]
Thus, $t$ must be equal to $0.56984...$ or $0.21508... \pm i \hskip 0.1em 1.30714...$ up to sign.
We set
\[
 t = 0.21508... - i \hskip 0.1em 1.30714...
\]
in the remaining.
Then $\xi \circ \rho_{\mathcal{A}}$ is in fact discrete and faithful.
Let $r_{1}, r_{2}, \cdots, r_{7} \in \Cpm$ denote the colors of regions with respect to $\mathcal{R}$ (see Figure \ref{fig:5_2}).
Suppose $r_1 = \begin{pmatrix} 1 \\ 0 \end{pmatrix}$.
According to the rule depicted in the right-hand side of Figure \ref{fig:coloring}, we have 
\[
 r_{1} = \begin{pmatrix} 1 \\ 0 \end{pmatrix}, \enskip r_{2} = \begin{pmatrix} 1 \\ t^{2} \end{pmatrix}, \enskip 
r_{3} = \begin{pmatrix} -t^{2} + 1 \\ t^{2} \end{pmatrix}, \enskip r_{4} = \begin{pmatrix} t^{4} - t^{2} + 1 \\ t^{4} \end{pmatrix},
\]
\[
 r_{5} = \begin{pmatrix} t^{4} + t^{2} + 1 \\ t^{4} \end{pmatrix}, \enskip r_{6} = \begin{pmatrix} -t^{6} - t^{4} + 1 \\ t^{2} (-t^{4} + 1) \end{pmatrix}, 
\enskip r_{7} = \begin{pmatrix} -t^{2} + 1 \\ t^{4} \end{pmatrix}.
\]

\begin{figure}[htb]
\begin{center}
\includegraphics[scale=0.32]{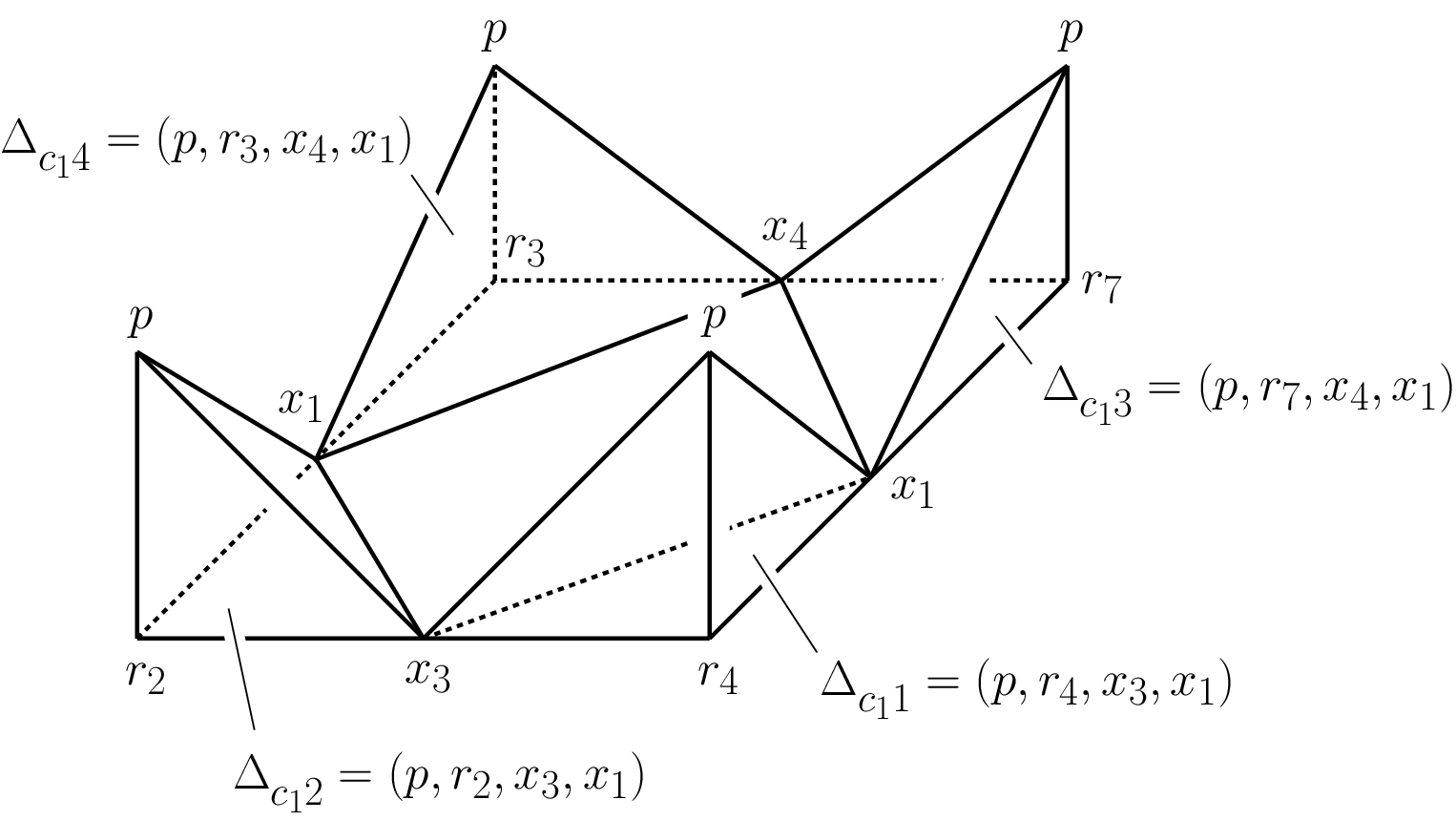}
\end{center}
\caption{Four simplices associated with $c_{1}$.}
\label{fig:crossing_1}
\end{figure}

We next compute $\langle [\mathrm{cvol}], [C(\mathcal{S})] \rangle$, that is the complex volume of $5_{2}$ knot.
Since $\langle [\mathrm{cvol}], [C(\mathcal{S})] \rangle = R(\widehat{\sigma}(\psi([C(\mathcal{S})])))$, we calculate $\widehat{\sigma}(\psi([C(\mathcal{S})])) \in \hatBC$.
We let $ p = \begin{pmatrix} 1 \\ -1 \end{pmatrix}$.
It is easy to see that $\varphi(C(\mathcal{S}))$ is not only an element of $C_{3}(\Cpm)_{G}$ 
but an element of $C_{3}^{h \neq}( \Cpm )_{G}$.
Thus, $\varphi(C(\mathcal{S}))$ is a representative of $\psi([C(\mathcal{S})])$.
For each crossing $c_{i}$, we have four simplices 
$\Delta_{c_{i} 1}, \Delta_{c_{i} 2}, \Delta_{c_{i} 3}, \Delta_{c_{i} 4} \in C_{3}^{h \neq}( \Cpm )_{G}$ 
derived from $\varphi(C(\mathcal{S}))$ (see Figure \ref{fig:crossing_1}).
Since $\widehat{\sigma}(\psi([C(\mathcal{S})]))$ is a signed sum of the combinatorial flattenings $\widehat{\sigma}(\Delta_{c_{i} j})$, 
it is sufficient to compute each $\widehat{\sigma}(\Delta_{c_{i} j})$.
For example, $\widehat{\sigma}(\Delta_{c_{1} 1})$ is calculated as follows.
By definition,
\begin{eqnarray*}
 w_{0} & = & \Log (\det(p, x_{1})) + \Log (\det(r_{4}, x_{3})) - \Log (\det(p, x_{3})) - \Log (\det(r_{4}, x_{1})) \\
       & = & \Log(1) + \Log(t^{3} - t) - \Log(-2t) - \Log(t^{4}) \\
       & = & - 0.816912... - i \hskip 0.1em 0.444187..., \\
 w_{1} & = & \Log (\det(p, x_{3})) + \Log (\det(r_{4}, x_{1})) - \Log (\det(p, r_{4})) - \Log (\det(a_{3}, a_{1})) \\
       & = & \Log(-2t) + \Log(t^{4}) - \Log(2t^{4} - t^{2} + 1) - \Log(-t) \\
       & = & - 0.344827... + i \hskip 0.1em 0.134887...
\end{eqnarray*}
Since the complex parameter $z$ of $\Delta_{c_{1} 1}$ is $[p : r_4 : x_3 : x_1] = \frac{t^{2} - 1}{2t^{4}}$, we have
\[
 \Log \left( z \right) = - 0.816912... + i \hskip 0.1em 2.69741..., \enskip - \Log \left( 1 - z \right) = - 0.344827... + i \hskip 0.1em 0.134887...
\]
Thus, $w_{0} = \Log \left( z \right) - \pi i$ and $w_{1} = - \Log \left( 1 - z \right)$.
Hence, $\widehat{\sigma}(\Delta_{c_{1} 1})$ is $[\frac{t^{2} - 1}{2t^{4}}, -1, 0]$.
By a straightforward calculation, we have
{\renewcommand{\arraystretch}{1.4}
\begin{center}
 \begin{tabular}[]{ll}
  $\widehat{\sigma}(\Delta_{c_{1} 1}) = [\frac{t^{2} - 1}{2t^{4}}; -1, 0]$, &
  $\widehat{\sigma}(\Delta_{c_{1} 2}) = [\frac{t^{2} - 1}{2t^{2}}; 0, 0]$, \\
  $\widehat{\sigma}(\Delta_{c_{1} 3}) = [\frac{t^{2} - 1}{t^{4}}; -1, 0]$, &
  $\widehat{\sigma}(\Delta_{c_{1} 4}) = [\frac{t^{2} - 1}{t^{2}}; 0, 1]$ \\
  $\widehat{\sigma}(\Delta_{c_{2} 1}) = [\frac{2t^{2} + 1}{t^{2}}; 0, -1]$, &
  $\widehat{\sigma}(\Delta_{c_{2} 2}) = [- \frac{2t^{2} + 1}{t^{4}}; 0, 0]$, \\
  $\widehat{\sigma}(\Delta_{c_{2} 3}) = [\frac{(2t^{2} + 1)(t^{2} + 1)^{2}}{t^{2}-1}; 0, -1]$, &
  $\widehat{\sigma}(\Delta_{c_{2} 4}) = [\frac{(2t^{2} + 1)(t^{6} + 2t^{4} - 2)}{t^{2}(t^{2} - 1)}; 0, -1]$, \\
  $\widehat{\sigma}(\Delta_{c_{3} 1}) = [\frac{2t^{4}}{t^{2} - 1}; 1, 0]$, &
  $\widehat{\sigma}(\Delta_{c_{3} 2}) = [\frac{2t^{4}}{t^{4} + t^{2} - 1}; 1, 0]$, \\
  $\widehat{\sigma}(\Delta_{c_{3} 3}) = [\frac{2t^{4}}{(2t^{2} + 1)(t^{2} - 1)}; 0, 0]$, &
  $\widehat{\sigma}(\Delta_{c_{3} 4}) = [\frac{2t^{4}}{(2t^{2} + 1)(t^{4} + t^{2} - 1)}; 0, -1]$, \\
  $\widehat{\sigma}(\Delta_{c_{4} 1}) = [\frac{t^{2}}{t^{2} - 1}; 0, 0]$, &
  $\widehat{\sigma}(\Delta_{c_{4} 2}) = [-\frac{t^{2}}{(t^{2} - 1)^{2}}; 0, 0]$, \\
  $\widehat{\sigma}(\Delta_{c_{4} 3}) = [-t^{2}; 0, -1]$, &
  $\widehat{\sigma}(\Delta_{c_{4} 4}) = [\frac{t^{2}}{t^{2} - 1}; 0, 0]$, \\
  $\widehat{\sigma}(\Delta_{c_{5} 1}) = [\frac{t^{2} - 1}{2t^{2} + 1}; 0, 0]$, &
  $\widehat{\sigma}(\Delta_{c_{5} 2}) = [\frac{t^{2} - 1}{(2t^{2} + 1)(t^{2} + 1)^{2}}; 0, 0]$, \\
  $\widehat{\sigma}(\Delta_{c_{5} 3}) = [- \frac{(t^{2} - 1)(t^{4} + t^{2} - 1)}{2t^{2}}; 0, 0]$, &
  $\widehat{\sigma}(\Delta_{c_{5} 4}) = [\frac{(t^{2} - 1)(t^{4} + t^{2} - 1)}{2t^{2}(t^{6} + 2t^{4} - 2)}; 0, 0]$.
 \end{tabular}
\end{center}
}

Evaluating a signed sum of these $\widehat{\sigma}(\Delta_{c_{i} j})$ with the map $R$, we have
\[
 \langle [\mathrm{cvol}], [C(\mathcal{S})] \rangle = R(\widehat{\sigma}(\psi([C(\mathcal{S})]))) = i \hskip 0.1em (2.82812... - i \hskip 0.1em 3.02412...).
\]
Therefore, the complex volume of $5_{2}$ knot is $i \hskip 0.1em (2.82812... - i \hskip 0.1em 3.02412...)$, as is known.


\begin{thebibliography}{99}
\bibitem{AG99}
N. Andruskiewitsch, M. Gra\~{n}a, 
{\it From racks to pointed Hopf algebras},
Adv. Math. 178 (2003), no. 2, 177--243.

\bibitem{CEGS99}
S. Carter, M. Elhamdadi, M. Gra\~{n}a, M. Saito,
{\it Cocycle knot invariants from quandle modules and generalized quandle homology},
Osaka J. Math. 42 (2005), no. 3, 499--541.

\bibitem{CES02}
J. S. Carter, M. Elhamdadi, M. Saito,
{\it Twisted quandle homology theory and cocycle knot invariants},
Algebr. Geom. Topol. 2 (2002), 95--135.

\bibitem{CJKLS03}
J. S. Carter, D. Jelsovsky, S. Kamada, L. Langford, M. Saito,
{\it Quandle cohomology and state-sum invariants of knotted curves and surfaces},
Trans. Amer. Math. Soc. 355 (2003), no. 10, 3947--3989.

\bibitem{CKS01}
J. S. Carter, S. Kamada, and M. Saito,
{\it Geometric interpretations of quandle homology},
J. Knot Theory Ramifications {\bf 10} (2001), no. 3, 345--386.

\bibitem{dupont}
J. L. Dupont,
{\it The dilogarithm as a characteristic class for flat bundles},
J. Pure Appl. Algebra 44 (1987), no. 1-3, 137--164.

\bibitem{dupont-zickert}
J. L. Dupont, C. Zickert,
{\it A dilogarithmic formula for the Cheeger-Chern-Simons class},
Geom. Topol. 10 (2006), 1347--1372. 

\bibitem{etingof-grana}
P. Etingof, M. Gra\~na,
{\it On rack cohomology},
J. Pure Appl. Algebra 177 (2003), no. 1, 49--59.

\bibitem{FRS95}
R.  Fenn, C.  Rourke and B. Sanderson, 
{\it Trunks and classifying spaces}, 
Appl. Categ. Structures 3 (1995), no. 4, 321--356.

\bibitem{inoue10_1}
A. Inoue,
{\it Quandle and hyperbolic volume},
Topology Appl. 157 (2010), no. 7, 1237--1245.

\bibitem{Inoue10_2}
A. Inoue,
{\it Knot quandles and infinite cyclic covering spaces},
Kodai Math. J. 33 (2010), no. 1, 116--122.

\bibitem{joyce}
D. Joyce, 
{\it A classifying invariant of knots, the knot quandle},
J. Pure Appl. Algebra 23 (1982), no. 1, 37--65.

\bibitem{kabaya}
Y. Kabaya,
{\it Cyclic branched coverings of knots and quandle homology}, 
Pacific J. Math. 259 (2012), no. 2, 315--347.

\bibitem{Kamada06}
S. Kamada,
{\it Quandles with good involutions, their homologies and knot invariants},
Intelligence of low dimensional topology 2006, 101--108, Ser. Knots Everything, 40, World Sci. Publ., Hackensack, NJ, 2007.

\bibitem{matveev}
S. Matveev,
{\it Distributive groupoids in knot theory} (Russian), 
Math. USSR-Sbornik 47 (1982), 73--83.

\bibitem{meyerhoff}
R. Meyerhoff, 
{\it Density of the Chern-Simons invariant for hyperbolic 3-manifolds}, 
Low-dimensional topology and Kleinian groups, 
London Math. Soc. Lecture Note Ser., 112 (1986), 217--239. 

\bibitem{neumann90}
W. Neumann,
{\it Combinatorics of triangulations and the Chern-Simons invariant for hyperbolic $3$-manifolds},
Topology '90 (Columbus, OH, 1990), 243--271.

\bibitem{neumann04}
W. Neumann,
{\it Extended Bloch group and the Cheeger-Chern-Simons class},
Geom. Topol. 8 (2004), 413--474.

\bibitem{neumann-yang}
W. Neumann, J. Yang, 
{\it Bloch invariants of hyperbolic $3$-manifolds}, 
Duke Math. J. 96 (1999), no. 1, 29--59.

\bibitem{niebrzydowski-przytycki}
M. Niebrzydowski, J. H. Przytycki,
{\it Homology operations on homology of quandles}, 
J. Algebra 324 (2010), no. 7, 1529--1548. 

\bibitem{weeks}
J. Weeks, 
{\it Computation of hyperbolic structures in knot theory}, 
Handbook of knot theory, 461--480, Elsevier B. V., Amsterdam, 2005.

\bibitem{zickert}
C. Zickert, 
{\it The volume and Chern-Simons invariant of a representation},  
Duke Math. J. 150 (2009), no. 3, 489--532. 

\end{thebibliography}
\end{document}